\documentclass[10pt]{amsart}
\usepackage{amsmath, amsthm, amssymb, euscript, enumerate}
\usepackage{graphicx, subfigure, epsfig}
\usepackage{mathtools, appendix, accents, color}
\usepackage[pdftex]{hyperref}
\usepackage{esvect}

\newcommand{\R}{{\mathbb R}}

\newcommand{\Z}{{\mathbb Z}}
\newcommand{\N}{{\mathbb N}}

\newcommand{\cS}{{\mathcal S}}

\newcommand{\e}{\varepsilon}

\newcommand{\norm}[1]{\left\Arrowvert {#1}\right\Arrowvert}
\newcommand{\inn}[2]{\left\langle{#1},{#2}\right\rangle}

\newcommand{\osc}{\operatornamewithlimits{osc}}

\newcommand{\diam}{\operatorname{diam}}
\newcommand{\dist}{\operatorname{dist}}

\newcommand{\tr}{\operatorname{tr}}

\theoremstyle{plain}
\newtheorem{theorem}{Theorem}[section]

\newtheorem{lemma}[theorem]{Lemma}
\newtheorem{proposition}[theorem]{Proposition}

\theoremstyle{definition}
\newtheorem{definition}[theorem]{Definition}

\theoremstyle{remark}

\newtheorem{remark}[theorem]{Remark}

\numberwithin{equation}{section}

\title[]{Uniform Estimates in Periodic Homogenization of Fully Nonlinear Elliptic Equations}
\author{Sunghan Kim}
\address{Department of Mathematics, KTH Royal Institute of Technology, 100 44 Stockholm, Sweden}
\email{sunghan@kth.se}
\author{Ki-Ahm Lee}
\address{Department of Mathematical Sciences, Seoul National University, Seoul 08826, Korea
\& Korea Institute for Advanced Study, Seoul 02455, Korea}
\email{kiahm@snu.ac.kr}
\thanks {S. Kim was supported by postdoctoral fellowship from Knut and Alice Wallenberg Foundation. K.-A. Lee was supported by NRF grant funded by the Korean government (MSIP) (NRF-2020R1A2C1A01006256). K.-A. Lee also holds a joint appointment with the Research Institute of Mathematics of Seoul National University.}

\begin{document}

\begin{abstract} 
This article is concerned with uniform $C^{1,\alpha}$ and $C^{1,1}$ estimates in periodic homogenization of fully nonlinear elliptic equations. The analysis is based on the compactness method, which involves linearization of the operator at each approximation step. Due to the nonlinearity of the equations, the linearized operators involve the Hessian of correctors, which appear in the previous step. The involvement of the Hessian of the correctors deteriorates the regularity of the linearized operator, and sometimes even changes its oscillating pattern. These issues are resolved with new approximation techniques, which yield a precise decomposition of the regular part and the irregular part of the homogenization process, along with a uniform control of the Hessian of the correctors in an intermediate level. The approximation techniques are even new in the context of linear equations. Our argument can be applied not only to concave operators, but also to certain class of non-concave operators. 
\end{abstract}


\maketitle

\tableofcontents


\section{Introduction}\label{section:intro}

This paper is concerned with the uniform estimates in periodic homogenization of fully nonlinear elliptic equations subject to a Dirichlet boundary condition, 
\begin{equation}\label{eq:main}
\begin{dcases}
F\left( D^2 u^\e, \frac{x}{\e}\right) = f & \text{in }\Omega,\\
u^\e  = g & \text{in } \partial\Omega.
\end{dcases}
\end{equation}
Here we establish uniform $C^{1,\alpha}$ and $C^{1,1}$ estimates for viscosity solutions $u^\e$ up to the boundary. The analysis is based on the compactness method \cite{AL1, AL2}, along with a new, sophisticated treatment on the nonlinear structure of the governing functional. 

Uniform regularity in the theory of homogenization has been of great interest for many years. A notable development was carried out by M. Avellaneda and F.-H. Lin in the series of works \cite{AL1,AL2}, where the compactness method is adopted in the framework of homogenization to prove uniform regularity estimates, which was primarily a technique for minimal surfaces and calculus of variations. Since then, the uniform estimates have been developed under various settings, such as singular integrals \cite{AL3}, Neumann boundary conditions \cite{KLS}, parabolic problems \cite{GS}, oscillating boundaries \cite{KP}, almost-periodic settings \cite{She2}. Some direct and constructive approaches were taken in \cite{CP, MS} and also in \cite{CL, LY} for homogenization problems of soft inclusions and highly oscillating obstacles. Recently, some new technique was developed for almost-periodic \cite{AS} and stationary ergodic problems \cite{AM} and Lipschitz domains \cite{She}, and it was also applied to nonlinear scalar equations \cite{WXZ} of divergence type. To the best of our knowledge, a sharp uniform interior estimate is available for nonlinear scalar equations of divergence type by \cite{WXZ}, but a sharp estimate up to the boundary is still left unsolved for the nonlinear equations.    

However, a sharp uniform estimate of $u^\e$ up to boundary for nonlinear problems has not yet been achieved in both divergence and non-divergence type equations, to the best of the authors' knowledge. Also the amount of the literature on non-divergence type equations is considerably smaller than that on the divergence type equations/systems. In this paper, we achieve a sharp uniform estimate, i.e., $C^{1,1}$ estimate up to the boundary for fully nonlinear elliptic problems of type \eqref{eq:main}. Our method can be carried out in various settings, such as nonlinear systems in divergence form, for which everywhere regularity is available, and parabolic problems. Nevertheless, we shall focus on the elliptic problems of non-divergence structure here, and leave the generalizations for the future. 

In what follows, we shall briefly explain our approach, and illustrate the main challenges arising from the nonlinearity of the governing operator $F$. 

The compactness method in a nutshell is a technique to approximate the solutions $u^\e$ by its limit profile, say $\bar u$, which is expected to have better regularity than $u^\e$, since $\bar u$ is meant to solve an effective problem that is homogeneous in small scales. One iteration step yields an approximation in an intermediate (yet universal) scale, say $\mu$, which could be much larger than the small scale, which is of order $\e$. Hence, one needs to iterate the approximation as many times as possible, e.g., $k$-times to reach $\mu^k \leq \e <\mu^{k-1}$. 

At each iteration step, we linearize the problem \eqref{eq:main} around the data obtained from the previous approximation step. Now if the governing operator $F$ were linear, i.e., $F(M,y) = a_{ij}(y)M_{ij}$, then the linearized operator coincides with the original one. Hence, one can only focus on the linearization of the source term and the boundary condition.

Nonetheless, if $F$ is nonlinear, the linearized operator records the data from the previous approximation step, and the data accumulates as the iteration continues. The accumulation effect itself is a generic issue in the analysis of the nonlinear equation; for instance, a similar issue appears in the Schauder estimate for fully nonlinear equations \cite[Theorem 8.1]{CC}. The difference here is that the accumulated data includes the (Hessian) of the interior and the boundary layer correctors, in order to bring back nice estimates for the limit profile $\bar u$ to $u^\e$. This is exactly where the analysis becomes difficult. 

To give a more precise picture of the accumulation effect, suppose that $\bar u$ satisfies a (universal) interior $C^{2,\bar\alpha}$ estimate at a given point. Then we can obtain along with a compactness argument that an interior $C^{2,\alpha}$ estimate, with $\alpha<\bar\alpha$, modulo an interior corrector, e.g.,
\begin{equation*}
\sup_{x\in B_\mu} \left| u^\e (x) - \frac{1}{2} \inn{M_1^\e x}{x} - \e^2 w_F\left( M_1^\e,\frac{x}{\e}\right) \right| \leq J^\e \mu^{2+\alpha},
\end{equation*} 
where $J^\e$ is the initial bound for $u^\e$ and $f$, $\mu$ is the universal intermediate scale in which the Hessian of $u^\e$ is approximated by $M_1^\e$, modulo the interior corrector $w_F(M_1^\e,\cdot)$ of $F$ at $M^\e$; more specifically, $w_F(M_1^\e,\cdot)$ is chosen by the periodic (viscosity) solution to the cell problem, 
\begin{equation*}
F(D^2 w + M_1^\e, y) = \bar F(M_1^\e)\quad \text{in }\R^n.
\end{equation*}
This suggests that in the next approximation step, we need to linearize the operator $F$ at $M_1^\e + D_y^2 w_F (M_1^\e,\cdot)$, instead of $M_1^\e$ alone, i.e., the linearized operator is given by 
\begin{equation}\label{eq:F1e}
F_1^\e (N,y) = \frac{F ( J^\e \mu^\alpha N + M_1^\e + D_y^2 w_F(M_1^\e,y),y) - F(M^\e + D_y^2 w_F(M_1^\e,y),y) }{J^\e \mu^\alpha}. 
\end{equation} 
Now we have to verify whether $F_1^\e$ satisfies all the structure conditions that were required in the compactness argument used in the first approximation step. This is the accumulation effect due to the nonlinear structure, since if $F$ were linear, we would have $F_1^\e = F$ which leaves us nothing to check.  

The uniform ellipticity and periodicity are preserved under the linearization. However, the regularity of the linearized operator $F_1^\e$ in the space variable now depends on that of $F$ in the matrix variable as well as the Hessian of the interior corrector $w_F(M_1^\e,\cdot)$ in the space variable. 

A quick remedy for this issue is to impose a strong assumption on $F$, such as $C^{1,1}$ regularity in the matrix variable and H\"older regularity in the space variable. Nevertheless, we present a new compactness argument (Lemma \ref{lemma:apprx}) that does not require any assumption other than uniform ellipticity, continuity and periodicity of the given operator; in particular, the argument does not require any control over the modulus of continuity. This argument is based on the observation (Lemma \ref{lemma:hom-stab}) that regardless of how irregular $F_1^\e$ is in its space variable, the effective functional $\bar F_1^\e$ stays homogeneous and has the same ellipticity bounds as those of ($F_1^\e$, hence of) $F$.

The strength of the general compactness argument (Lemma \ref{lemma:apprx}) is that we obtain a uniform interior $C^{1,1}$ estimate for $u^\e$ up to a deleted neighborhood of size $\e$ (Theorem \ref{theorem:int-C11} (i)), for any uniformly elliptic, periodic functional $F$ that admits $C^2$-regular correctors. Moreover, the modulus of continuity of the Hessian of the correctors is not involved in the estimate. This provides us a better understanding of the homogenization process, as the regularity estimate decomposes $u^\e$ into regular part and irregular part. In particular, the irregular part is regularized under a stronger assumption that the governing operator $F$ is H\"older continuous in its space variable (Theorem \ref{theorem:int-C11} (ii)). In other words, a (standard) uniform interior $C^{1,1}$ estimate of $u^\e$ is established, with more regularity assumption on $F$. Let us remark that such a decomposition of $u^\e$ is even new in the context of linear equations. 

A more challenging issue appears in the analysis of the sharp, uniform boundary estimate. Again we linearize the problem around the approximating data from the previous iteration step, and the data accumulates due to the nonlinearity of the given operator $F$. This times, however, the accumulated data does not only include the interior correctors but also the boundary layer correctors. The problem is that the boundary layer correctors do not adhere the same periodic pattern in the rapid oscillation as the interior correctors do. Thus, once the boundary layer correctors appear in the linearized operator, the operator no longer oscillates periodically in the space variable. This implies a qualitative change of the oscillating nature in the next approximation step. 

The authors encountered a similar issue in \cite{KL2}, while studying the higher-order convergence rates in periodic homogenization of oscillating initial data. Here the problem is more delicate because the oscillatory patterns of the boundary layer correctors are more ambiguous than those of the initial layer correctors considered in \cite{KL2}. Still, we shall need an analogous treatment that the effect from the boundary layer correctors dissipates fast as we stay away from the boundary layer, but at the same time it can be controlled uniformly up to the boundary with respect to certain norm, which is less sharper, yet sufficient for the analysis.

To explain the issue regarding the boundary layer correctors in more details, let us suppose that we have approximated $u^\e$ by its limit profile $\bar u$, which satisfies a boundary $C^{2,\bar\alpha}$ estimate at $0\in \partial\Omega$, as 
\begin{equation*}
\sup_{\Omega_\mu} \left| u^\e(x) - \frac{1}{2} \inn{M_1^\e x}{x} - \e^2 w_F \left(M_1^\e,\frac{x}{\e}\right) - \zeta_1^\e(x)\right| \leq J^\e \mu^{2+\alpha}, 
\end{equation*}
with some $\alpha<\bar\alpha$, where $J^\e$ is the initial bound for $u^\e$, $f$ and $g$, $\Omega_\mu = \Omega\cap B_\mu$ and $\zeta_1^\e$ is the boundary layer corrector. Here the boundary layer corrector $\zeta_1^\e$ is given by the (viscosity) solution to 
\begin{equation*}
\begin{dcases}
F \left( M_1^\e + D_y^2 w_F\left(M_1^\e, \frac{x}{\e}\right) + D_x^2 \zeta_1^\e,\frac{x}{\e} \right) = \bar F(M_1^\e) & \text{in }\Omega_1,\\
\zeta_1^\e = -\e^2 w_F\left(M_1^\e, \frac{x}{\e}\right) & \text{on }\partial\Omega_1.
\end{dcases}
\end{equation*} 
In other words, $\zeta_1^\e$ is chosen so as to correct the error on the boundary layer, $\Gamma_1 = \partial\Omega\cap B_1$, left from the interior correction by $\e^2 w_F(M_1^\e,\frac{\cdot}{\e})$ in $\Omega_1$. Then the linearized operator for the next approximation step will be of the form,
\begin{equation*}
\begin{split} 
G_1^\e(N,x,y) &= \frac{1}{J^\e \mu^\alpha} F( J^\e\mu^\alpha N + M_1^\e + D_y^2 w_F(M_1^\e, y) + D_x^2 \zeta_1^\e (\mu x), y) \\
& \quad - \frac{1}{J^\e\mu^\alpha} F ( M_1^\e + D_y^2 w_F(M_1^\e,y) + D_x^2 \zeta_1^\e (\mu x),y).
\end{split} 
\end{equation*} 

It is noteworthy that $G_1^\e$ depends on $x$, and the dependence is through $D^2 \zeta_1^\e$, which is an irregular term. If $F$ were linear in the matrix variable, say $F(M,y) = a_{ij}(y) M_{ij}$, then the effect of $D^2 \zeta_1^\e$ is cancelled out, and we end up with $G_1^\e (N,x,y) = a_{ij}(y) N_{ij} = F(N,y)$. For this reason, the main focus for the linear problems is set to the uniform control of the linearized boundary condition, as shown in \cite{AL2}.

If the dependence of $G_1^\e$ on $x$ were regular and uniform (e.g., of class $C^\alpha$ uniformly), this would not be an issue either. However, the dependence here is far from being uniform continuous, since both the interior equation and the boundary condition for $\zeta_1^\e$ are rapidly oscillating in $\frac{x}{\e}$. In fact, we do not have a uniform $L^\infty$ bound on $D^2 \zeta_1^\e$ up to the boundary at this stage, as it is what we are aiming for. This suggests that we need to control $D^2 \zeta_1^\e$ with some weaker uniform estimates. 

We overcome this difficulty (in Lemma \ref{lemma:apprx-bdry}) as follows. First, we observe that the effect of $D^2 \zeta_1^\e$ becomes negligible in the interior, or more precisely, 
\begin{equation*}
|D^2\zeta_1^\e(x)| = O\left( \frac{\e^2}{\dist(x,\Gamma_1)^2}\right) \quad\text{for }x\in\Omega_{1/2}, 
\end{equation*}
where $\Omega_{1/2} = \Omega\cap B_{1/2}$. This is proved with the uniform interior $C^{1,1}$ estimate (Theorem \ref{theorem:int-C11}) that we establish prior to the boundary analysis, together with an elementary {\it a priori} estimate $\zeta_1^\e = O(\e^2)$ in $\Omega_1$, which follows immediately from the boundary condition for $\zeta_1^\e$ and the maximum principle.

On the other hand, the boundary condition for $\zeta_1^\e$ also implies that $D_x^2 \zeta_1^\e = - D_y^2 w_F(M_1^\e,\frac{x}{\e}) = O(1)$ on $\Gamma_1$, from which we observe a uniform control of $D^2 \zeta_1^\e$ up to the boundary in $L^p$ sense, for any $p>n$ large, i.e.,
\begin{equation*}
\norm{D^2 \zeta_1^\e}_{L^p(\Omega_{1/2})} = O(1).
\end{equation*}
This is based on a uniform $W^{2,p}$ estimate up to the boundary (Proposition \ref{proposition:bdry-W2p}), which we derive by combining the above uniform interior $C^{1,1}$ estimate with the uniform boundary $C^{1,\beta}$ estimate (Theorem \ref{theorem:bdry-C1a}) in a standard way. 

These two observations allow us to decompose the aperiodic, linearized operator $G_1^\e$ into a sum of a periodic operator and a source term, which is uniformly bounded in the $L^p$ space and dissipates away from the boundary layer in the $L^\infty$ sense. For this reason, we can maintain the periodic nature of the original problem \eqref{eq:main} during the entire iteration scheme (Lemma \ref{lemma:iter-bdry}), through a suitable compactness argument (Lemma \ref{lemma:apprx-bdry}). 

As a final remark, let us discuss about the effective operators in the homogenization of fully nonlinear problems. It is by now standard \cite{E} that if $F$ is uniformly elliptic and periodic, then the problem \eqref{eq:main} is homogenized into 
\begin{equation*}
\begin{cases}
\bar F (D^2 \bar u) = f &\text{in }\Omega,\\
\bar u = g &\text{on } \partial\Omega,
\end{cases}
\end{equation*}  
for some homogeneous, elliptic, or the so-called effective operator $\bar F$, in the sense that $u^\e \to \bar u$ uniformly over $\bar \Omega$. Hence, if $\bar u$ is a regular solution, e.g., $\bar u \in C^{2,\alpha}$, then one may expect to use this information to derive better regularity for $u^\e$, which constitutes the basic idea of the compactness method.

Unlike linear problems, the homogeneous equation $H (D^2 u ) = 0$ does not always admit classical solutions. According to Nadirashvili and Vladut \cite{NV}, there are homogeneous, elliptic functionals $H$ that admit viscosity solutions belonging to $C^{1,\beta}\setminus C^{1,1}$ for some $\beta\in(0,1)$. For this reason, we shall focus ourselves on the periodically oscillating operator $F$ such that not only the homogeneous equation $F(D^2 u ,y_0) = 0$ obtained by ``freezing coefficients'' at each $y_0$, but also the effective one $\bar F(D^2 u) = 0$ admit classical, or more precisely, $C^{2,\alpha}$ solutions; this notion will be made precise in Definition \ref{definition:basic} and Definition \ref{definition:class}. 

One may ask under which condition on $F$, viscosity solutions to the homogeneous equation $\bar F ( D^2 u ) = 0$ verify a universal interior $C^{2,\alpha}$ estimate. Let us stress that the class of such operators $F$ is non-void. A typical example would be concave operators $F$ (i.e., $F(M,y)$ is concave in the matrix variable $M$ and periodic in the space variable $y$). Note that if $F$ is concave, then so is the effective operator $\bar F$, according to \cite{E}; then by the Evans-Krylov theory \cite{CC}, interior $C^{2,\alpha}$ estimates are available for both homogeneous equations $F(D^2 v,y_0) = 0$ and $\bar F (D^2 u ) = 0$. 

However, there has not yet been any observation of a class beyond concave operators in the literature, for which the effective operators inherit the Evans-Krylov type estimates. Here, we present certain class (Proposition \ref{proposition:ex}) of non-concave, periodic operators $F$ whose effective operators $\bar F$ admit interior $C^{2,\alpha}$ estimates. This shows that the uniform regularity estimates we establish in this paper are applicable for a wide class of periodic functionals $F$.

On a different note, we believe that it deserves independent interests to find the largest class of such periodic operators. This can be reformulated as follows: if $F$ is a uniformly elliptic, periodic functional such that $F(\cdot,y_0)$ admits an interior $C^{2,\alpha}$ estimate independent of $y_0\in\R^n$, does the effective functional $\bar F$ also admit an interior $C^{2,\alpha'}$ estimate, possibly with some other $\alpha'$? In an abstract level, it is a question concerning the structures of periodically oscillating functionals that are preserved during the homogenization process. Uniform ellipticity and concavity are some typical properties that are preserved, as shown in \cite{E}. More recently, the authors have proved in the series of papers \cite{KL1,KL2,KL3} that higher regularity of $F$ in the matrix variable is also inherited to $\bar F$, i.e., if $F(\cdot,y_0)\in C_{loc}^{k,1}$ uniformly in $y_0\in\R^n$, then $\bar F\in C_{loc}^{k-1,1}$. From this point of view, here we only present a partial answer to the question, and aim to give a more complete picture in the future. 

We are in a position to state the main theorem of this paper. Here the classes $S_0$, $S_1$, $S_2$, $R_0$, $R_1$ regarding the structure and regularity conditions of governing functionals are defined in Definition \ref{definition:class}, and $\bar{F}$ is the effective functional as in Definition \ref{definition:eff}.

\begin{theorem}\label{theorem:main}
Let $F\in S_0(\lambda,\Lambda)$ be a periodic functional on $\cS^n\times\R^n$, $\Omega\subset\R^n$ be a bounded domain with $\partial\Omega\in C^2$, $f\in L^\infty(\Omega)$ and $g\in C^{1,\alpha}(\Omega_1)$, for some $0<\lambda\leq\Lambda$ and $0<\alpha<1$. Let $u^\e$ be the viscosity solution of \eqref{eq:main} for each $\e>0$. Let $\bar\kappa$, $\bar\gamma$, $\kappa$, $\gamma$, $\bar{c}$ and $\bar\alpha$ be additional positive parameters with $\gamma<\bar\gamma\leq 1$ and $\alpha<\bar\alpha\leq 1$. 
\begin{enumerate}[(i)]
\item
If $F\in S_1(\lambda,\Lambda,\bar\kappa,\bar\gamma)\cap R_0(\kappa,\gamma)$ having $\bar{F} \in S_1(\lambda,\Lambda,\bar{c},\bar\alpha)$, then $u^\e \in C^{1,\alpha_*}(\Omega)$, with $\alpha_* = \min\{\alpha,\bar\gamma\}$ and 
\begin{equation}\label{eq:main-C1a}
\norm{u^\e}_{C^{1,\alpha_*}(\Omega)} \leq C_1 \left( \norm{f}_{L^\infty(\Omega)} + \norm{g}_{C^{1,\alpha}(\partial\Omega)} \right),
\end{equation}
where $C_1$ depends only on $n$, $\lambda$, $\Lambda$, $\bar\kappa$, $\bar\gamma$, $\kappa$, $\gamma$, $\bar{c}$, $\bar\alpha$, $\alpha$, $\diam(\Omega)$ and the maximal curvature of $\partial\Omega$. 
\item
If $F\in S_2(\lambda,\Lambda,\bar\kappa,\bar\gamma)\cap R_1(\kappa,\gamma)$ having $\bar{F}\in S_2(\lambda,\Lambda,\bar{c},\bar\alpha)$, $\partial\Omega \in C^{2,\bar\alpha}$, $f\in C^\alpha(\Omega)$ and $g\in C^{2,\alpha}(\partial \Omega)$, then $u^\e \in C^{1,1}(\Omega)$ and 
\begin{equation}\label{eq:main-C11}
\norm{u^\e}_{C^{1,1}(\Omega)} \leq C_2,
\end{equation}
with $C_2$ depending only on $n$, $\lambda$, $\Lambda$, $\bar\kappa$, $\bar\gamma$, $\kappa$, $\gamma$, $\bar{c}$, $\bar\alpha$, $\alpha$, $\norm{f}_{C^\alpha(\Omega)}$, $\norm{g}_{C^{2,\alpha}(\partial\Omega)}$, $\diam(\Omega)$ and the $C^{2,\bar\alpha}$ character of $\partial\Omega$. 
\end{enumerate}
\end{theorem}

It should be remarked that the assumption of Theorem \ref{theorem:main} (i) on the class $S_1$ is always satisfied, if the exponents $\bar\gamma$ and $\bar\alpha$ are small enough, since a homogeneous, elliptic functional always admits an interior $C^{1,\delta}$ estimate for some $\delta\in(0,1)$, depending only on the space dimension and the ellipticity constants \cite[Corollary 5.7]{CC}. On the other hand, the assumption of Theorem \ref{theorem:main} (ii) on the class $S_2$, especially regarding $\bar{F}$, is essential, since a nonlinear functional fails to have an interior $C^{2,\delta}$ estimate in general.

Let us briefly introduce the notation used throughout this article. Number $n\geq 1$ will always denote the space dimension. $\cS^n$ is the space of all symmetric $n\times n$ matrices. For $N\in\cS^n$, by $N\geq 0$ we shall indicate that all the eigenvalues of $N$ are nonnegative. By $|N|$ we shall denote the $(L^2,L^2)$ norm of $N$, that is, $|N| = (\sum_{i,j=1}^n N_{ij})^{1/2}$. Given a set $A\subset\R^n$ and a point $x\in\R^n$, $d(x,A) = \min(\dist(x,A),1)$. By a domain, we refer to an open connected set. For a domain $\Omega$, by $\partial\Omega$ we denote the topological boundary of $\Omega$. Occasionally, we shall also write by $\Gamma$ the boundary $\partial\Omega$. We shall also write $\Omega_r(z) = \Omega\cap B_r(z)$, $\Gamma_r(z) = \Gamma\cap B_r(z)$, $\Omega_r = \Omega_r(0)$ and $\Gamma_r = \Gamma_r(0)$. Note that $\partial\Omega_r(z) = \Gamma_r(z) \cup (\Omega\cap \partial B_r(z))$. Given a domain $\Omega$, we define H\"older spaces, $C^{k,\alpha}(\Omega)$ and $C_{loc}^{k,\alpha}(\Omega)$ in the usual sense. Also when $k=0$, we shall simply denote them by $C^\alpha(\Omega)$ and respectively $C_{loc}^\alpha(\Omega)$. 

The paper is organized as follows. In Section \ref{section:prelim}, we list up some preliminaries required for this paper. In Section \ref{section:stab} we observe new stability results that will be used in the compactness method in the subsequent sections. Section \ref{section:int-C1a} is devoted to the uniform interior $C^{1,\alpha}$ estimates (Theorem \ref{theorem:int-C1a}). Section \ref{section:bdry-C1a} is concerned with the uniform interior $C^{1,\alpha}$ estimates (Theorem \ref{theorem:bdry-C1a}). In Section \ref{section:int-C11}, we study the uniform interior $C^{1,1}$ estimates (Theorem \ref{theorem:int-C11}), and in Section \ref{section:bdry-C11}, we consider the uniform boundary $C^{1,1}$ estimates (Theorem \ref{theorem:bdry-C11}). The proof of Theorem \ref{theorem:main} is omitted, as it follows straightforwardly from these four theorems. Finally, in Section \ref{section:ex}, we present certain class of non-concave functionals satisfying the assumptions of Theorem \ref{theorem:main}. 


\section{Preliminaries}\label{section:prelim}

Let us first introduce some terminologies that will appear throughout this paper. 

\begin{definition}\label{definition:basic}
Let $Y\subset\R^n$ be a domain, and $F$ be a continuous functional on $\cS^n\times Y$. 
\begin{itemize}
\item $F$ is said to be uniformly elliptic with constants $\lambda$ and $\Lambda$, if 
\begin{equation}\label{eq:F-ellip}
\lambda |N| \leq F(M+N,y) - F(M,y) \leq \Lambda |N|,\quad N\in\cS^n,N\geq 0, 
\end{equation}
for any $M\in\cS^n$ and $y\in Y$. 
\item When $Y = \R^n$, $F$ is said to be periodic, if  
\begin{equation}\label{eq:F-peri}
F(M, y +k) = F(M,y),\quad k\in\Z^n,
\end{equation}
for any $M\in\cS^n$ and $y\in\R^n$. 
\item $F$ is said to have zero source term, if 
\begin{equation}\label{eq:F-0}
F(0,y) = 0,\quad y\in Y. 
\end{equation}
\item $F$ is said to have an interior $C^{k,\gamma}$ estimate with constant $\kappa$, if for any $y_0\in Y$ with $a = F(0,y_0)$ and any viscosity solution $v\in C(\bar{B}_1)$ of $F(D^2 v,y_0) = a$ in $B_1$, one has $v\in C^{k,\gamma}(B_{1/2})$ and 
\begin{equation*}
\norm{v}_{C^{k,\gamma}(B_{1/2})} \leq \kappa \norm{v}_{L^\infty(B_1)}. 
\end{equation*}
\end{itemize}
\end{definition}

We shall consider several classes of functionals defined as follows.

\begin{definition}\label{definition:class}
Let $F$ be a continuous functional on $\cS^n\times Y$, and let $\lambda$, $\Lambda$, $\bar\kappa$, $\bar\gamma$, $\kappa$ and $\gamma$, with $\lambda\leq\Lambda$ and $\gamma<\bar\gamma\leq 1$. Class $S_0$, $S_1$ and $S_2$, regarding structure conditions of governing functionals, are defined as follows. 
\begin{itemize}
\item $F\in S_0(\lambda,\Lambda)$, if $F$ is uniformly elliptic with ellipticity constants $\lambda$ and $\Lambda$, and has zero source term. 
\item $F \in S_1(\lambda,\Lambda,\bar\kappa,\bar\gamma)$, if $F\in S_0(\lambda,\Lambda)$ and has an interior $C^{1,\bar\gamma}$ estimate with constant $\bar\kappa$. 
\item $F\in S_2(\lambda,\Lambda,\bar\kappa,\bar\gamma)$, if $F\in S_0(\lambda,\Lambda)$ and $(M,y) \mapsto (F(M +N ,y) - F(N,y))$ has an interior $C^{2,\bar\gamma}$ estimate with constant $\bar\kappa$, uniformly for all $N\in\cS^n$. 
\end{itemize}
Classes $R_0$, $R_1$ and $R_1'$ will refer to regularity conditions of governing functionals, as follows.
\begin{itemize}
\item $F\in R_0(\kappa,\gamma)$, if $F\in C(\cS^n\times Y)$ and 
\begin{equation}\label{eq:F-Ca}
|F(M,y_1) - F(M,y_2)| \leq \kappa |M| |y_1-y_2|^{\gamma},
\end{equation}
for any $M\in\cS^n$ and $y_1,y_2\in Y$. 
\item $F\in R_1(\kappa,\gamma)$, if $F\in R_0(\kappa,\gamma)$ and $F(\cdot,y)\in C^1(\cS^n)$ for any $y\in Y$ satisfying 
\begin{equation}\label{eq:F-C11}
\left| D_M F (M_1,y_1) - D_M F (M_2,y_2) \right| \leq \kappa ( |M_1 - M_2| + |y_1-y_2|^\gamma),
\end{equation}
for any pair $(M_i,y_i)\in\cS^n\times Y$ with $i \in\{1,2\}$. 
\end{itemize}
If the parameters are not important, or are well understood from the context, we shall omit them and simply write $S_0$, $S_1$, and so on.
\end{definition}

Next we define the class of $C^{2,\alpha}$ domains.

\begin{definition}\label{definition:domain}
Given a domain $\Omega$ with $0<\bar\sigma\leq 1$ and $\bar\tau>0$, we shall say $\Omega\in D(\bar\tau,\bar\sigma)$, if $\partial\Omega$ is locally a $C^{2,\bar\tau}$-graph around the origin, with $C^{2,\bar\tau}$-character being controlled by $\bar\sigma$. More specifically, we indicate that $0\in\partial\Omega$ and there are a function $\phi\in C^{2,\bar\sigma}(\R^{n-1})$ with
\begin{equation*}
\phi(0) = |D_T \phi(0)| = 0,\quad \norm{\phi}_{C^{2,\bar\sigma}(\R^{n-1})} \leq \bar\tau,
\end{equation*}
and a rotation $\Phi:\R^n\to\R^n$, which maps the inward unit normal $\nu$ to $\partial\Omega$ at the origin to $e_n$, such that with $\Omega_1 = \Omega\cap B_1$,
\begin{equation*}
\Phi(\Omega_1) \subset \{(x',x_n)\in B_1: x_n > \phi(x')\}.
\end{equation*}
Let us remark that as $\Phi$ being a rotation, $|D_T\phi(0)| = 0$ does not imply that the inward normal vector $\nu$ to $\partial\Omega$ at the origin is the same with $e_n$. 
\end{definition}

Let us introduce the notation for effective functional and interior corrector, whose existence and uniqueness of effective functionals are well understood, and we refer to \cite{E} for a proof. 

\begin{definition}\label{definition:eff}
Let $F$ be a uniformly elliptic, continuous and periodic functional on $\cS^n\times\R^n$. By $\bar{F}$ we shall always denote the functional on $\cS^n$ satisfying for each $M\in\cS^n$, $\bar{F}(M)$ is the unique number for which there exists a periodic viscosity solution to
\begin{equation}\label{eq:cell}
F(M + D_y^2 w,y) = \bar{F}(M)\quad\text{in }\R^n.
\end{equation}
$\bar{F}$ will be called the effective functional of $F$. Equation \eqref{eq:cell} will be called the cell problem associated with $F$. 

Moreover, by $w_F$ we shall denote a functional on $\cS^n\times\R^n$ such that for each $M\in\cS^n$, $w_F(M,\cdot)$ is the unique periodic viscosity solution to \eqref{eq:cell} satisfying
\begin{equation}\label{eq:w-0}
w_F(M,0) = 0. 
\end{equation}
$w_F$ will be called the normalized interior corrector associated with $F$. 
\end{definition}

In what follows, we shall list up some basic facts that will be used frequently throughout the paper. Since the proofs are elementary, and we believe are well understood by the experts, we shall not present them here. 

First, we state the closeness of the classes $S_1$ and $S_2$.

\begin{lemma}\label{lemma:limit} Fix $i\in\{1,2\}$. Let $\{F_k\}_{k=1}^\infty\subset S_i(\lambda,\Lambda,\bar\kappa,\bar\gamma)$ be a sequence of continuous functionals on $\cS^n\times Y$, with $Y\subset\R^n$ a domain, and suppose that $F_k\to F$ locally uniformly on $\cS^n\times Y$, as $k\to\infty$, for some functional $F$. Then $F\in S_i(\lambda,\Lambda,\bar\kappa,\bar\gamma)$. 
\end{lemma}

This ensures the well-definedness of normalized interior correctors. Let us list up some basic properties that will be used throughout this article. 

\begin{lemma}\label{lemma:Fb}
Let $F\in S_0(\lambda,\Lambda)$ be a continuous and periodic functional on $\cS^n\times\R^n$. 
\begin{enumerate}[(i)]
\item $\bar{F}\in S_0(\lambda,\Lambda)$ on $\cS^n$, in the sense that $\bar{F}(0) = 0$ and 
\begin{equation}\label{eq:Fb-ellip}
\lambda |M| \leq \bar{F}( M + N ) - \bar{F}(N) \leq \Lambda |M|,\quad N\in\cS^n, N \geq 0,
\end{equation}
for any $M\in\cS^n$. 
\item There are $0<\alpha<1$ and $C_0>0$, depending only on $n$, $\lambda$ and $\Lambda$, such that for each $M\in\cS^n$, $w_F(M,\cdot)\in C^\alpha(\R^n)$ and 
\begin{equation}\label{eq:wF-Ca}
\norm{w_F(M,\cdot)}_{C^\alpha(\R^n)} \leq C_0 |M|.
\end{equation}
\item Moreover, if $F\in R_0(\kappa,\gamma)$, then there are an exponent $0<\beta<1$, depending only on $n$, $\lambda$ and $\Lambda$, and a constant $C_1>0$, depending only on $n$, $\lambda$, $\Lambda$, $\kappa$ and $\gamma$, such that for each $M\in\cS^n$, $w_F(M,\cdot) \in C^{1,\beta}(\R^n)$ and 
\begin{equation}\label{eq:wF-C1a}
\norm{w_F(M,\cdot)}_{C^{1,\beta}(\R^n)} \leq C_1 |M|. 
\end{equation}
\item Assume further that $F\in S_2(\lambda,\Lambda,\bar\kappa,\bar\gamma)\cap R_0(\kappa,\gamma)$. Then there is $C_2>0$, depending only on $n$, $\lambda$, $\Lambda$, $\bar\kappa$, $\bar\gamma$, $\kappa$ and $\gamma$, such that for each $M\in\cS^n$, $w_F(M,\cdot) \in C^{2,\gamma}(\R^n)$ and 
\begin{equation}\label{eq:wF-C2a}
\norm{ w_F (M,\cdot)}_{C^{2,\gamma}(\R^n)} \leq C_2|M|.
\end{equation} 
\end{enumerate}
\end{lemma}

In the subsequent two lemmas, we shall define, given a periodic functional $F \in S_2\cap R_0$ on $\cS^n\times\R^n$, 
\begin{equation}\label{eq:FMm}
\begin{split}
F_{M,\mu} (N,y) &= \frac{F( \mu N + M + D_y^2 w_F(M,y),y) - F( M + D_y^2 w_F(M,y),y)}{\mu} \\
& = \frac{F( \mu N + M + D_y^2 w_F(M,y),y) - \bar{F}(M)}{\mu}.
\end{split}
\end{equation}
First we observe how such a scaling affects the regularity of the associated functional. 

\begin{lemma}\label{lemma:GF} 
Let $F\in S_2(\lambda,\Lambda,\bar\kappa,\bar\gamma)\cap R_0(\kappa,\gamma)$ be a periodic functional on $\cS^n\times\R^n$. Then for each $M\in\cS^n$ and $\mu>0$, $F_{M,\mu} \in S_2(\lambda,\Lambda,\bar\kappa,\bar\gamma)$. Moreover, 
\begin{equation*}
F_{M,\mu}\in R_0\left(\kappa\left(1+\frac{C_0|M|}{\mu}\right),\gamma\right),
\end{equation*}
where $C_0$ depends only on $n$, $\lambda$, $\Lambda$, $\bar\kappa$, $\bar\gamma$, $\kappa$ and $\gamma$. Moreover, if $F\in R_1(\kappa,\gamma)$, then for any $M\in\cS^n$ and $\mu>0$,
\begin{equation*}
F_{M,\mu} \in R_0 (\kappa (C_0|M| + 1 ),\gamma).
\end{equation*}
\end{lemma}

Next lemma amounts to the change of the effective functional and the normalized corrector under the given scaling. 

\begin{lemma}\label{lemma:eff-cor} 
Let $F\in S_0(\lambda,\Lambda)$ be a periodic functional on $\cS^n\times\R^n$, and suppose that $w_F(M,\cdot) \in C^2(\R^n)$ for any $M\in\cS^n$. Then for each $M\in\cS^n$ and $\mu>0$, one has 
\begin{equation*}
\bar{F}_{M,\mu} (N) = \frac{\bar{F}(\mu N + M) - \bar{F}(M)}{\mu},\quad\text{and}
\end{equation*}
\begin{equation*}
w_{F_{M,\mu}} (N,y) = \frac{w_F(\mu N + M,y) - w_F(M,y)}{\mu},
\end{equation*}
for all $N\in\cS^n$ and $y\in\R^n$. 
\end{lemma}


\section{Stability in Homogenization}\label{section:stab}

Let us present a stability result of viscosity solutions in periodic homogenization problem. A key difference from the classical stability theory such as \cite[Proposition 4.10]{CC} is that we do not assume the family of governing functionals to be convergent. Although the proof involves a minor modification of the classical argument \cite[Theorem 3.3]{E}, we find the assertion itself interesting, and for the sake of completeness, we shall contain the full arguments here. 

\begin{lemma}\label{lemma:hom-stab}
Let $\{F_k\}_{k=1}^\infty\subset S_0(\lambda,\Lambda)$ be a sequence of continuous periodic functional on $\cS^n\times\R^n$. Let $\Omega$ be a domain, $\{f_k\}_{k=1}^\infty$ a sequence of continuous functions on $\Omega$ converging locally uniformly to $f$. Given a sequence $\{\e_k\}_{k=1}^\infty$ of positive numbers decreasing to zero, suppose that $\{u_k\}_{k=1}^\infty$ is a uniformly bounded sequence of viscosity solutions to 
\begin{equation*}
F_k \left( D^2 u_k,\frac{x}{\e_k} \right) = f_k\quad\text{in }\Omega.
\end{equation*} 
Then there are a functional $\bar{F}\in S_0(\lambda,\Lambda)$ on $\cS^n$, and a function $\bar{u}\in C(\Omega)$ such that $\bar{F}_k \to \bar{F}$ locally uniformly on $\cS^n$ and $u_k\to \bar{u}$ locally uniformly in $\Omega$, after extracting a subsequence if necessary. Moreover, $\bar{u}$ is a viscosity solution of 
\begin{equation*}
\bar{F}(D^2 \bar{u}) = f \quad\text{in }\Omega.
\end{equation*}
\end{lemma}

\begin{proof}
Since $F_k\in S_0(\lambda,\Lambda)$ on $\cS^n\times\R^n$, the effective functional $\bar{F}_k$ on $\cS^n$ also belongs to $S_0(\lambda,\Lambda)$, c.f. Lemma \ref{lemma:Fb} (i). In particular, it ensures the uniform Lipschitz continuity of the sequence $\{\bar{F}_k\}_{k=1}^\infty$ on $\cS^n$. Hence, after extracting a subsequence, if necessary, there is a functional $\bar{F}$ on $\cS^n$ to which $\{\bar{F}_k\}_{k=1}^\infty$ converges locally uniformly on $\cS^n$. It is clear that $\bar{F}$ is a functional on $\cS^n$ belonging to $S_0(\lambda,\Lambda)$.

Since $f_k$ is locally bounded in $\Omega$ and $u_k$ is uniformly bounded in $\Omega$, it follows from the Krylov-Safanov theory \cite[Proposition 4.9]{CC} that $\{u_k\}_{k=1}^\infty$ is a bounded sequence in $C^\alpha(K)$, for each compact $K\subset\Omega$, where $0<\alpha<1$ depends only on $n$, $\lambda$ and $\Lambda$. Here we used that the ellipticity constants of $F_k$ are fixed by $\lambda$ and $\Lambda$. Therefore, after extracting a subsequence if necessary, $u_k\to \bar{u}$ locally uniformly in $\Omega$. 

Hence, we are left with proving the main assertion of this lemma that $\bar{F}(D^2 \bar{u}) = f$ in $\Omega$ in the viscosity sense. The rest of the proof follows closely the argument of \cite[Theorem 3.3]{E}. 

Without loss of generality, let us assume that $\bar{F}_k \to \bar{F}$ and $u_k\to \bar{u}$ along the full sequence as $k\to\infty$. In order to prove that $\bar{u}$ is a viscosity subsolution, we assume towards a contradiction that for some $x_0\in\Omega$, there exists $\phi\in C^2(\Omega)$ touching $\bar{u}$ strictly from above at $x_0$, but
\begin{equation*}
\bar{F}(D^2 \phi(x_0)) \leq f(x_0) - 3\theta,
\end{equation*}
for some $\theta>0$. Since $\bar{F}_k (D^2 \phi(x_0)) \to \bar{F}(D^2 \phi(x_0))$, and $f_k(x_0) \to f(x_0)$ as $k\to\infty$, one has 
\begin{equation}\label{eq:false2-stab}
\bar{F}_k(D^2 \phi(x_0)) \leq f_k(x_0) - 2\theta,
\end{equation}
for all sufficiently large $k$'s. 

As $\bar{F}_k$ being the effective functional of $F_k$, with the normalized interior corrector $w_{F_k}$, one has 
\begin{equation}\label{eq:cell-stab}
F_k ( D_x^2\phi(x_0) + D_y^2 w_{F_k}(D_x^2 \phi(x_0),y),y) = \bar{F}_k(D_x^2 \phi(x_0))\leq f_k(x_0) - 3\theta
\end{equation}
in the viscosity sense in $\R^n$. Since $F_k \in S_0(\lambda,\Lambda)$, we know from Lemma \ref{lemma:Fb} (ii) that
\begin{equation}\label{eq:w-stab}
\norm{w_{F_k} (D_x^2 \phi(x_0),\cdot)}_{L^\infty(\R^n)} \leq C_0(1 + | D_x^2 \phi(x_0)| ),
\end{equation}
for all $k=1,2,\cdots$, where $C_0$ depends only on $n$, $\lambda$ and $\Lambda$. 

Define 
\begin{equation*}
\phi_k(x) = \phi(x) + w_{F_k}\left(D_x^2 \phi(x_0), \frac{x}{\e}\right),
\end{equation*}
and let us claim that $\phi_k$ is a viscosity solution to
\begin{equation}\label{eq:claim-stab}
F_k\left( D^2 \phi_k,\frac{x}{\e_k}\right) \leq f_k(x) - \theta\quad\text{in }B_r(x_0),
\end{equation}
provided that $r$ is sufficiently small. 

Let us choose $r$ small enough such that $\overline{B_r(x_0)} \subset \Omega$, and 
\begin{equation}\label{eq:phi-stab}
\Lambda | D^2 \phi(x) - D^2 \phi(x_0) | + |f_k(x) - f_k(x_0)| \leq \theta \quad\text{ for any }x\in B_r(x_0),
\end{equation}
by utilizing the continuity of $D^2 \phi$ in $\Omega$ and the assumption that $f_k\to f$ locally uniformly in $\Omega$. Moreover, taking $r$ smaller if necessary, we may assume that
\begin{equation}\label{eq:phi2-stab}
\bar{u}(x_0) = \phi(x_0), \quad \max_{\partial B_r(x_0)} (\bar{u} - \phi) < 0,
\end{equation}
since $\phi$ was taken to be a function touching $\bar{u}$ strictly from above at $x_0$. 

In order to prove the claim, let $x_1\in B_r(x_0)$ be arbitrary, and suppose that $\psi\in C^2(B_r(x_0))$ is a function touching $\phi_k$ from below at $x_1$. Then the auxiliary function 
\begin{equation*}
\eta(y) = \frac{1}{\e_k^2} (\psi (\e_k y ) - \phi (\e_k y)) = \frac{1}{\e_k^2} (\psi(\e_k y) - \phi_k (\e_k y)) + w_{F_k}\left( D_x^2 \phi(x_0), \frac{x}{\e_k}\right)
\end{equation*}
belongs to the class $C^2(B_{\e_k^{-1} r}( \e_k^{-1}x_0))$ and touches $w_k$ from below at $y_1 = \e_k^{-1} x_1$. As $w_{F_k} (D_x^2 \phi(x_0),\cdot)$ being a viscosity solution to \eqref{eq:cell-stab}, we have
\begin{equation*}
F_k ( D_x^2 \phi(x_0) + D_y^2 \eta(y_1), y_1) \leq f_k(x_0) - 2\theta,
\end{equation*}
and since $D_y^2 \eta(y_1) = D_x^2 \psi(x_1) - D_x^2 \phi(x_1)$, one may deduce, from \eqref{eq:phi-stab} and the fact that $\Lambda$ is the upper ellipticity bound of $F_k$, that
\begin{equation*}
\begin{split}
F_k \left( D_x^2 \psi (x_1),\frac{x_1}{\e_k}\right) &\leq F_k (D_x^2 \phi(x_0) + D_y^2 \eta(y_1), y_1) + \Lambda |D_x^2 \phi(x_1) - D_x^2 \phi(x_0)|\\
& \leq f_k(x_0) - |f_k(x_1) - f_k(x_0)| - \theta \\
& \leq f_k(x_1) - \theta,
\end{split}
\end{equation*}
which proves the claim \eqref{eq:claim-stab}, for any large $k$'s. 

Now that $F_k$ is a uniformly elliptic functional, the comparison principle yields that 
\begin{equation*}
u_k(x_0) - \phi_k(x_0) \leq \max_{\partial B_r(x_0)} (u_k - \phi_k),
\end{equation*}
for all sufficiently large $k$'s. Passing to the limit with $k\to\infty$, and using the uniform convergence of $u_k\to \bar{u}$ on $\bar{B}_r(x_0)$ together with the uniform estimate \eqref{eq:w-stab}, one arrives at 
\begin{equation*}
\bar{u}(x_0) - \phi(x_0) \leq \max_{\partial B_r(x_0)} (\bar{u} - \phi)
\end{equation*}
which is a contradiction to \eqref{eq:phi2-stab}. Hence, one should have $\bar{F}(D^2 \bar{u}) \geq f$ in $\Omega$ in the viscosity sense. The reverse inequality can also be proved in a similar argument, and we omit the details. 
\end{proof}

Next let us state a stability result in order to formulate a suitable approximation lemma for boundary estimate. Now we require the sequence $\{F_k\}_{k=1}^\infty$ of functionals to lie in the class $S_2\cap R_0$, with fixed parameters, since we are going to allow certain aperiodic perturbation proportional to the size of the Hessian variable.

\begin{lemma}\label{lemma:hom-stab2}
Let $\{f_k\}_{k=1}^\infty$, $\{\e_k\}_{k=1}^\infty$ and $\Omega$ be as in Lemma \ref{lemma:hom-stab}. Let $\{F_k\}_{k=1}^\infty\subset S_2(\lambda,\Lambda,\bar\kappa,\bar\gamma)\cap R_0(\kappa,\gamma)$ be a sequence of periodic functionals on $\cS^n\times\R^n$, and suppose that $\{G_k\}_{k=1}^\infty\subset S_0(\lambda,\Lambda)$ is another sequence of periodic functionals on $\cS^n\times\Omega\times \R^n$ such that for any compact $K\subset\Omega$ and $\eta>0$, there exists a sufficiently large $k_0\geq 1$ for which 
\begin{equation*}
| G_k (N,x,y) - F_k (N,y) | \leq \eta |N|,\quad \text{if }k\geq k_0,
\end{equation*}
for all $N\in\cS^n$, $x\in K$ and $y\in\R^n$. Under this assumption, let $\{u_k\}_{k=1}^\infty$ be a sequence of locally uniformly bounded viscosity solutions to 
\begin{equation*}
G_k \left( D^2 u_k ,x,\frac{x}{\e_k}\right) = f_k \quad \text{in }\Omega. 
\end{equation*} 
Then there are a functional $F \in S_2(\lambda,\Lambda,\bar\kappa,\bar\gamma)\cap R_0(\kappa,\gamma)$ on $\cS^n\times\R^n$ and a function $\bar{u} \in C(\Omega)$ such that $F_k \to F$ locally uniformly in $\cS^n\times\R^n$ and $u_k \to \bar{u}$ locally uniformly in $\Omega$, after extracting a subsequence, and that
\begin{equation*}
\bar{F}(D^2 \bar{u}) = f \quad \text{in }\Omega.
\end{equation*} 
\end{lemma}

\begin{proof}
Note from Lemma \ref{lemma:limit} that the class $S_2(\lambda,\Lambda,\bar\kappa,\bar\gamma)$ is closed under a locally uniform convergence, while the uniform estimate \eqref{eq:F-Ca} on the H\"older continuity of $F_k$ yields the compactness of the sequence $\{F_k\}_{k=1}^\infty$. Thus, we obtain $F\in S_2(\lambda,\Lambda,\bar\kappa,\bar\gamma)\cap R_0(\kappa,\gamma)$ such that $F_k \to F$ locally uniformly in $\cS^n\times\R^n$, after extracting a subsequence. Then it is clear that $\bar{F}_k \to \bar{F}$ locally uniformly on $\cS^n$, along a subsequence, where $\bar{F}$ now is the effective functional corresponding to $F$. In addition, one may follow the same argument in Lemma \ref{lemma:hom-stab} and obtain a function $\bar{u}\in C(\Omega)$ to which $\{\bar{u}_k\}_{k=1}^\infty$ converges locally uniformly in $\Omega$, after extracting a subsequence. Here we shall assume without loss of generality that all the convergence listed above holds along the full sequence. 

Assume to the contrary that there exist some $x_0\in\Omega$ and certain $\phi\in C^2(\Omega)$ such that $\phi$ touches $\bar{u}$ strictly from above at $x_0$ but 
\begin{equation*}
\bar{F}(D^2 \phi(x_0)) = f(x_0) - 4\theta,
\end{equation*}
for some $\theta>0$, and hence 
\begin{equation}\label{eq:hom-stab-bdry1}
\bar{F}_k (D^2 \phi(x_0)) \leq f(x_0) - 3\theta,
\end{equation}
for all sufficiently large $k$'s. Let us denote by $N$ the matrix $D^2 \phi(x_0)$, so that $w_{F_k}(N,\cdot)$ is a periodic viscosity solution of 
\begin{equation}\label{eq:hom-stab-bdry2}
F_k (N + D_y^2 w_{F_k}(N,y), y) = \bar{F}(N) \leq f(x_0) - 3\theta\quad\text{in }\R^n.
\end{equation}
From the assumption that $F_k \in S_2(\lambda,\Lambda,\bar\kappa,\bar\gamma)\cap R_0(\kappa,\gamma)$, Lemma \ref{lemma:Fb} (iv) yields that $w_{F_k} (N,\cdot) \in C^{2,\gamma}(\R^n)$ and satisfies 
\begin{equation*}
\norm{D_y^2 w_{F_k}(N,\cdot)}_{L^\infty(\R^n)} \leq C_0 |N|,
\end{equation*} 
where $C_0$ depends only on $n$, $\lambda$, $\Lambda$, $\bar\kappa$, $\bar\gamma$, $\kappa$ and $\gamma$. Thus, the auxiliary function 
\begin{equation*}
\phi_k (x) = \phi(x) + \e_k^2 w_{F_k} \left(N, \frac{x}{\e_k}\right)
\end{equation*}
belongs to $C^2(\Omega)$, and moreover
\begin{equation}\label{eq:hom-stab-bdry3}
\norm{D^2 \phi_k}_{L^\infty(K)} \leq \norm{D^2 \phi}_{L^\infty(K)} + C_0 |N| \leq (1+C_0)\norm{D^2 \phi}_{L^\infty(K)}, 
\end{equation}
for any compact set $K\subset\Omega$ containing $x_0$; note that the last inequality follows from $|N| = |D^2 \phi(x_0)| \leq \norm{D^2 \phi}_{L^\infty(K)}$. 

Let $k_0\geq 1$ be such that $\{f_k\}_{k=k_0}^\infty$ is uniformly continuous over $K$. Such a number $k_0$ exists, since $f_k \in C(\Omega)$ and  $f_k \to f$ locally uniformly in $\Omega$ as $k\to\infty$. Now from the assumption that $\phi$ touches $\bar{u}$ strictly from above, the hypothesis that $\phi\in C^2(\Omega)$, and the choice of $k_0$, we can choose $r>0$ small enough such that
\begin{equation}\label{eq:hom-stab-bdry7}
\min_{\partial B_r(x_0)} (\phi - \bar{u}) > 0,
\end{equation} 
and 
\begin{equation}\label{eq:hom-stab-bdry4}
\Lambda \norm{D^2 \phi - N}_{L^\infty(B_r(x_0))} + \sup_{k\geq k_0} \norm{f_k  - f_k(x_0)}_{L^\infty(B_r(x_0))} \leq \theta.
\end{equation} 
Next we choose $\eta$ so small that
\begin{equation}\label{eq:hom-stab-bdry5}
(1+C_0) \eta \norm{D^2 \phi}_{L^\infty(B_r(x_0))} \leq \theta,
\end{equation} 
with $C_0$ as in \eqref{eq:hom-stab-bdry3}, and then correspondingly select a sufficiently large $k_1\geq k_0$ such that 
\begin{equation}\label{eq:hom-stab-bdry6}
|G_k( N, x, y ) - F_k(N,y)| \leq \eta |N|, \quad\text{if } k\geq k_0, 
\end{equation}
for all $N\in\cS^n$, $x\in B_r(x_0)$ and $y\in\R^n$. Then for any $k\geq k_0$, it follows from \eqref{eq:hom-stab-bdry4}, \eqref{eq:hom-stab-bdry5}, \eqref{eq:hom-stab-bdry6}, \eqref{eq:hom-stab-bdry3} and \eqref{eq:hom-stab-bdry2} that 
\begin{equation*}
\begin{split}
G_k\left( D^2 \phi_k (x), x,\frac{x}{\e_k}\right) &\leq F_k \left( D^2 \phi_k(x),\frac{x}{\e_k}\right) + \eta \norm{D^2 \phi_k}_{L^\infty(B_r(x_0))} \\
&\leq F_k \left( N + D_y^2 w_{F_k} \left(N,\frac{x}{\e_k}\right),\frac{x}{\e_k}\right) + \Lambda |D^2 \phi(x)- N| + \theta \\
&\leq f(x_0) - \theta - |f_k(x) -f_k(x_0)| \\
&\leq f(x) - \theta, 
\end{split}
\end{equation*} 
for any $x\in B_r(x_0)$, in the classical sense. Consequently, one may follow the argument in Lemma \ref{lemma:hom-stab} and arrive at a contradiction to \eqref{eq:hom-stab-bdry7}. This finishes the proof. 
\end{proof}


\section{Interior $C^{1,\alpha}$ Estimate}\label{section:int-C1a}

In this section, we study uniform interior $C^{1,\alpha}$ estimates of viscosity solutions to fully nonlinear elliptic equations. 
The main result is stated as follows. 

\begin{theorem}\label{theorem:int-C1a}
Let $F\in S_0(\lambda,\Lambda)$ be a periodic functional on $\cS^n\times\R^n$ with $\bar{F} \in S_1(\lambda,\Lambda,\bar{c},\bar\alpha)$, $f\in L^\infty(B_1)$ and let $u^\e$ be a viscosity solution to 
\begin{equation}\label{eq:main-int-C1a-pde}
F \left( D^2 u^\e,\frac{x}{\e}\right) = f\quad\text{in }B_1. 
\end{equation}
Let $\alpha\in(0,\bar\alpha)$ be arbitrary. 
\begin{enumerate}[(i)]
\item Then for any $x_0\in B_{1/2}$, and any $\e\in(0,\frac{1}{8})$, there exists an affine function $l_{x_0}^\e$, with $l_{x_0}^\e (x_0) = u^\e(x_0)$, such that 
\begin{equation}\label{eq:int-C1a-almst}
\begin{split}
|\nabla l_{x_0}^\e| + \sup_{x\in B_{3/4}(x_0)\setminus \overline{B_\e(x_0)}} \frac{| u^\e (x) - l_{x_0}^\e (x)|}{|x-x_0|^{1+\alpha}} \leq C_1 \left( \norm{u^\e}_{L^\infty(B_1)} + \norm{f}_{L^\infty(B_1)} \right), 
\end{split}
\end{equation}
where $C_1>0$ depends only on $n$, $\lambda$, $\Lambda$, $\bar{c}$, $\bar\alpha$ and $\alpha$. 
\item If $F \in  S_1(\lambda,\Lambda,\bar\kappa,\bar\gamma)\cap R_0(\kappa,\gamma)$, then $u^\e \in C^{1,\min(\alpha,\gamma)}(B_{1/2})$ and 
\begin{equation}\label{eq:int-C1a}
\norm{u^\e}_{C^{1,\min(\alpha,\gamma)}(B_{1/2})} \leq C_2 \left(\norm{u^\e}_{L^\infty(B_1)} + \norm{f}_{L^\infty(B_1)}\right),
\end{equation}
where $C_2>0$ depends in addition on $\kappa$, $\gamma$, $\bar\kappa$ and $\bar\gamma$.
\end{enumerate}
\end{theorem}

The first part of the theorem shows that without any control on the modulus of continuity for $F$ in its spatial variable, we can still obtain a $C^{1,\alpha}$ estimate all the way up to an $\e$-neighborhood, essentially proving that the irregular part of $u^\e$ due to the rapid oscillation in the governing operator can be isolated in $\e$-cubes. Such an observation was deduced from the stability result established in Lemma \ref{lemma:hom-stab}) that allows the compactness method to work in a very general framework. Now if one assume stronger regularity on $F$, we can also control the oscillating behavior of $u^\e$ inside the $\e$-cube, and recover a full $C^{1,\alpha}$ estimate (part (ii) above). 

Let us also mention that the almost $C^{1,\alpha}$ estimate (Theorem \ref{theorem:int-C1a} (i)) will be used essentially in the subsequent analysis, where we establish uniform $C^{1,1}$ and $W^{2,p}$ estimates. 

Our analysis begins with an approximation lemma. 

\begin{lemma}\label{lemma:apprx-C1a}
Let $\bar\alpha\in(0,1]$ and $\alpha \in (0,\bar\alpha)$ be arbitrary. One can choose small positive constants $\mu$ and $\bar\e$, depending only on $n$, $\lambda$, $\Lambda$, $\bar{c}$, $\bar\alpha$ and $\alpha$, such that if $\e\leq \bar\e$, $F\in S_0(\lambda,\Lambda)$ is a periodic functional on $\cS^n\times\R^n$ with $\bar{F} \in S_1(\lambda,\Lambda,\bar{c},\bar\alpha)$, and $f\in L^\infty(B_1)$, $u^\e\in C(B_1)$ satisfying
\begin{equation}\label{eq:f-Linf}
\norm{f}_{L^\infty(B_1)} \leq \e,
\end{equation}
\begin{equation}\label{eq:ue-pde-C1a}
F\left( D^2u^\e, \frac{x}{\e}\right) =f\quad\text{in }B_1,
\end{equation}
\begin{equation}\label{eq:ue-Linf-C1a}
u^\e(0) = 0 \quad\text{and}\quad \norm{u^\e}_{L^\infty(B_1)} \leq 1, 
\end{equation}
then 
\begin{equation}\label{eq:apprx-C1a}
\sup_{x\in B_\mu} \left| u^\e (x) - \frac{1}{\omega_n\mu^n}\int_{\partial B_\mu} u(z) \frac{\inn{x}{z}}{|z|} \,d\sigma \right| \leq \mu^{1+\alpha},
\end{equation}
where $\omega_n$ is the volume of the unit ball in $\R^n$. 
\end{lemma}

\begin{proof} 
The proof is the essentially same with the argument in \cite[Lemma 4]{AL2}, whence we shall provide the details for the points where the arguments differ. Suppose by way of contradiction that this lemma fails to hold. Then there are a sequence $\{\e_k\}_{k=1}^\infty$ of positive numbers decreasing to zero such that for each $k=1,2,\cdots$, a sequence $\{F_k\}_{k=1}^\infty\subset S_0(\lambda,\Lambda)$ of periodic functionals on $\cS^n\times\R^n$ having $\bar{F}\in S_1(\lambda,\Lambda,\bar{c},\bar\alpha)$ on $\cS^n$, sequences $\{f_k\}_{k=1}^\infty\subset L^\infty(B_1)$ and $\{u_k\}_{k=1}^\infty\subset C(\overline{B_1})$ satisfying
\begin{equation}\label{eq:fk-Linf-false}
\norm{f_k}_{L^\infty(B_1)} \leq \e_k,
\end{equation}
\begin{equation}\label{eq:uk-pde-C1a-false}
F_k \left( D^2 u_k, \frac{x}{\e_k} \right) = f_k\quad\text{in }B_1,
\end{equation}
\begin{equation}\label{eq:uk-Linf-C1a-false}
u_k(0) = 0\quad\text{and}\quad \norm{u_k}_{L^\infty(B_1)} \leq 1,
\end{equation}
such that 
\begin{equation}\label{eq:apprx-C1a-false}
\sup_{x\in B_\mu} \left| u_k(x) - \frac{1}{\omega_n\mu^n}\int_{\partial B_\mu} u_k(z) \frac{\inn{x}{z}}{|z|} \,d\sigma \right| > \mu^{1+\alpha}.
\end{equation}

Since $F_k$ has the same ellipticity constants for all $k$, and $u_k$ is a uniformly bounded solution of \eqref{eq:uk-pde-C1a-false}, it follows from the Krylov-Safanov theory \cite[Proposition 4.9]{CC} that $\{u_k\}_{k=1}^\infty$ is uniformly bounded in $C_{loc}^\beta(B_1)$, for some $0<\beta<1$, depending only on $n$, $\lambda$ and $\Lambda$. Without loss of generality, one may assume that there is $\bar{u} \in C_{loc}^\beta$ for which 
\begin{equation}\label{eq:uk-ub-Ca-false}
u_k \to \bar{u} \quad\text{locally uniformly in }B_1. 
\end{equation} 
In particular, it follows from the uniform bound \eqref{eq:uk-Linf-C1a-false} of $u_k$ that
\begin{equation}\label{eq:ub-Linf-C1a-false}
\norm{\bar{u}}_{L^\infty(B_1)} \leq 1. 
\end{equation}

On the other hand, \eqref{eq:fk-Linf-false} implies that $f_k\to 0$ in $L^\infty(B_1)$ as $k\to\infty$. Hence, Lemma \ref{lemma:hom-stab} yields a functional $\bar{F}\in S_0(\lambda,\Lambda)$ on $\cS^n$, which can be characterized as a limit of $\{\bar{F}_k\}_{k=1}^\infty$ under the locally uniform convergence on $\cS^n$, such that 
\begin{equation}\label{eq:ub-pde-C1a-false}
\bar{F}(D^2 \bar{u}) = 0 \quad\text{in }B_1,
\end{equation}
in the viscosity sense. Due to the assumption that $\bar{F}_k \in S_1(\lambda,\Lambda,\bar{c},\bar\alpha)$ on $\cS^n$, Lemma \ref{lemma:limit} ensures that $\bar{F} \in S_1(\lambda,\Lambda,\bar{c},\bar\alpha)$ as well. In view of \eqref{eq:ub-Linf-C1a-false} and \eqref{eq:ub-pde-C1a-false}, we have  
\begin{equation}\label{eq:ub-C1a-false}
\norm{\bar{u}}_{C^{1,\bar\alpha}(B_{3/4})} \leq \bar{c}.
\end{equation} 
From this point the proof is the same with \cite[Lemma 4]{AL2}, so we skip the details.
\end{proof} 

Next we proceed with an iteration lemma.

\begin{lemma}\label{lemma:iter-C1a}
Let $\bar\alpha$, $\bar{c}$, $\alpha$, $\mu$ and $\bar\e$ be as in Lemma \ref{lemma:apprx-C1a}. If $\e\leq \mu^{k-1}\bar\e$ for some integer $k\geq 1$, $F\in S_0(\lambda,\Lambda)$ is a periodic function on $\cS^n\times\R^n$ having $\bar{F} \in S_1(\lambda,\Lambda,\bar{c},\bar\alpha)$ on $\cS^n$, $f\in L^\infty(B_1)$ and $u^\e$ is a viscosity solution of
\begin{equation}\label{eq:ue-pde-C1a-gen}
F\left( D^2 u^\e,\frac{x}{\e} \right) = f\quad\text{in }B_1,
\end{equation}
then for each $l\in\{1,\cdots,k\}$, there exists $a_l^\e\in \R^n$, with $a_0^\e = 0$, such that 
\begin{equation}\label{eq:ake-C1a}
|a_l^\e - a_{l-1}^\e| \leq C_0 \mu^{(l-1)\alpha} J^\e,\quad\text{and}
\end{equation}
\begin{equation}\label{eq:iter-C1a}
\sup_{x\in B_{\mu^l}} | u^\e (x) - u^\e (0) - \inn{a_l^\e}{ x} | \leq J^\e \mu^{l(1+\alpha)},
\end{equation}
where $C_0$ is a positive constant, depending only on $n$, $\lambda$, $\Lambda$, $\bar{c}$, $\bar\alpha$ and $\alpha$, 
\begin{equation}\label{eq:Je-C1a}
J^\e = \norm{u^\e}_{L^\infty(B_1)} + \frac{1}{\bar\e} \norm{f}_{L^\infty(B_1)}.
\end{equation}
\end{lemma}

\begin{proof}
Let $C_0$ be a constant to be determined later. The case $\e\leq \bar\e$ (i.e., $k=1$ in the statement) can be treated easily from Lemma \ref{lemma:apprx-C1a}, by applying the lemma to $(J^\e)^{-1} u^\e$. Hence, we shall assume that Lemma \ref{lemma:iter-C1a} is true for the case when $\e\leq \mu^{k-1}\bar\e$ for some $k\geq 2$, and we shall attempt to prove that it continues to hold when $\e\leq \mu^k\bar\e$ for the same integer $k$. In particular, from the induction hypothesis, we already have $a_l^\e\in\R^n$, for $l\in\{1,\cdots,k\}$, such that both \eqref{eq:ake-C1a} and \eqref{eq:iter-C1a} are satisfied. 

Consider the scaled functions, 
\begin{equation}\label{eq:uke-fke-C1a}
u_k^\e (x) = \frac{u^\e(\mu^k x) - u^\e(0) - \mu^k \inn{a_k^\e}{x}}{J^\e \mu^{k(1+\alpha)}},\quad f_k^\e(x) = \frac{\mu^{k(1-\alpha)}f(\mu^k x)}{J^\e},
\end{equation}
which are defined for $x\in B_1$. From the induction hypothesis \eqref{eq:iter-C1a} for $u^\e$, which is a viscosity solution of \eqref{eq:ue-pde-C1a-gen}, we know that
\begin{equation}\label{eq:uke-pde-C1a}
F_k^\e \left( D^2 u_k^\e,\frac{\mu^k x}{\e} \right) = f_k^\e\quad\text{in }B_1,
\end{equation}
\begin{equation}\label{eq:uke-Linf-C1a}
u_k^\e (0) = 0 \quad\text{and}\quad \norm{u_k^\e}_{L^\infty(B_1)}\leq 1, 
\end{equation}
where we wrote 
\begin{equation}\label{eq:Fe-C1a}
F_k^\e (N,y) = \frac{\mu^{k(1-\alpha)}}{J^\e} F \left( \frac{J^\e}{\mu^{k(1-\alpha)}} N, y\right).
\end{equation}
It is clear that $F_k^\e$ is a periodic functional on $\cS^n\times\R^n$ belonging to $S_0(\lambda,\Lambda)$. Moreover, Lemma \ref{lemma:eff-cor} ensures that 
\begin{equation*}
\bar{F}_k^\e (N) = \frac{\mu^{k(1-\alpha)}}{J^\e} \bar{F} \left( \frac{J^\e}{\mu^{k(1-\alpha)}}N \right),
\end{equation*}
and hence it is clear that $\bar{F}_k^\e \in S_1(\lambda,\Lambda,\bar{c},\bar\alpha)$ as well on $\cS^n$. 

On the other hand, since $0<\alpha<1$, it follows that 
\begin{equation}\label{eq:fke-Linf-C1a}
\norm{f_k^\e}_{L^\infty(B_1)} \leq \frac{1}{J^\e} \norm{f}_{L^\infty(B_{\mu^k})}\leq \bar\e. 
\end{equation} 
Thus, one can apply Lemma \ref{lemma:apprx-C1a} and obtain 
\begin{equation}\label{eq:iter-C1a-1}
\sup_{x\in B_\mu} \left| u_k^\e (x) - \frac{1}{\omega_n\mu^n} \int_{\partial B_\mu} u_k^\e (z) \frac{\inn{x}{z}}{|z|} \,d\sigma \right| \leq \mu^{1+\alpha}
\end{equation}
Rewriting this inequality in terms of $u^\e$, we see that
\begin{equation}\label{eq:iter-C1a-2}
\sup_{x\in B_{\mu^{k+1}}} |u^\e (x) - u^\e (0) - a_{k+1}^\e \cdot x | \leq J^\e \mu^{(k+1)(1+\alpha)},
\end{equation}
with
\begin{equation}\label{eq:ake-re2}
a_{k+1}^\e = a_k^\e + \frac{J^\e \mu^{k\alpha}}{\omega_n\mu^n} \int_{\partial B_\mu} u_k^\e(z) \frac{z}{|z|}\,d\sigma.
\end{equation}
From \eqref{eq:iter-C1a} we deduce that 
\begin{equation*}
|a_{k+1}^\e - a_k^\e | \leq n J^\e\mu^{k\alpha - 1},
\end{equation*}
verifying \eqref{eq:ake-C1a} for $k+1$. Thus, the proof is finished with $C_0 = n\mu^{-1}$, which certainly depends only on $n$, $\lambda$, $\Lambda$, $\bar{c}$, $\bar\alpha$ and $\alpha$. 
\end{proof}

We are ready to prove the uniform interior $C^{1,\alpha}$ estimates.

\begin{proof}[Proof of Theorem \ref{theorem:int-C1a}]
The first assertion of Theorem \ref{theorem:int-C1a} follows directly from Lemma \ref{lemma:iter-C1a}. The second part of the statement is also similar with the proof of \cite[Theorem 1 (i)]{AL2}, due to the same lemma. Hence we omit the details. 
\end{proof}


\section{Boundary $C^{1,\alpha}$ Estimate}\label{section:bdry-C1a}

The main objective of this section is to prove uniform boundary $C^{1,\alpha}$ estimates. 

\begin{theorem}\label{theorem:bdry-C1a}
Let $F\in S_1(\lambda,\Lambda,\bar\kappa,\bar\gamma)\cap R_0(\kappa,\gamma)$ be a periodic functional on $\cS^n\times\R^n$, $\Omega$ be domain with $\Gamma\in C^2$ containing the origin, $f\in L^\infty(\Omega)$ and $g\in C^{1,\alpha}(\Gamma)$, for some $0<\alpha\leq 1$. Suppose that $u^\e$ is a viscosity solution of 
\begin{equation}\label{eq:main-bdry-C1a-pde}
\begin{dcases}
F \left( D^2 u^\e, \frac{x}{\e} \right) = f & \text{in }\Omega_1,\\
u^\e = g & \text{on }\Gamma_1.
\end{dcases}
\end{equation}
Then $u^\e \in C^{1,\min(\alpha,\gamma)}(\Omega_{1/2})$ and 
\begin{equation}\label{eq:bdry-C1a}
\norm{u^\e}_{C^{1,\min(\alpha,\gamma)}(\Omega_{1/2})} \leq C  \left(  \norm{u^\e}_{L^\infty(\Omega_1)} + \norm{f}_{L^\infty(\Omega_1)} + \norm{g}_{C^{1,\bar\alpha}(\Gamma_1)} \right),
\end{equation}
where $C$ depends only on $n$, $\lambda$, $\Lambda$, $\bar\kappa$, $\bar\gamma$, $\kappa$, $\gamma$, $\alpha$ and the maximal curvature of $\Gamma_1$. 
\end{theorem}

It is noteworthy that the boundary $C^{1,\alpha}$ estimate is rather straightforward than the interior case. For instance, it follows from the work \cite{SS} by L. Silvestre and B. Sirakov that $u^\e$ satisfies a uniform boundary $C^{1,\alpha}$ estimate with $\alpha$ depending only on $n$, $\lambda$ and $\Lambda$, since the equation \eqref{eq:main-bdry-C1a-pde} essentially belongs to the Pucci class. This is also the reason why we do not formulate the discrete boundary $C^{1,\alpha}$ estimate in contrast with the interior case. 

What is new here is that one can improve the H\"older exponent as much as one desires. This is due to the fact that homogeneous functionals always have interior $C^{2,\bar\alpha}$ estimate for some universal $\bar\alpha$, again due to \cite{SS}. 

As in the previous subsection, we shall begin with an approximation lemma. 

\begin{lemma}\label{lemma:apprx-bdry-C1a}
Let $0<\alpha<1$ be arbitrary, and $\Omega$ be a domain such that $\Gamma\in C^2$ contains the origin and the maximal curvature of $\Gamma_1$ is bounded by $1$. One can choose $0<\mu,\bar\e\leq\frac{1}{2}$, depending only on $n$, $\lambda$, $\Lambda$ and $\alpha$, such that if $F\in S_0(\lambda,\Lambda)$ is a periodic functional on $\cS^n\times\R^n$, $f\in L^\infty(\Omega_1)$, $g\in C^{1,\alpha}(\Gamma_1)$, $u^\e\in C(\overline{\Omega_1})$ satisfying 
\begin{equation}\label{eq:f-Linf-bdry-C1a}
\norm{f}_{L^\infty(\Omega_1)} \leq \bar\e,
\end{equation}
\begin{equation}\label{eq:g-C1a-bdry-C1a}
g(0) = |D_T g(0)| = 0, \quad \norm{g}_{C^{1,\alpha}(\Gamma_1)} \leq \bar\e,
\end{equation}
\begin{equation}\label{eq:ue-pde-bdry-C1a-apprx}
\begin{dcases}
F \left( D^2 u^\e ,\frac{x}{\e}\right) = f,\quad |u^\e|\leq 1 & \text{in }\Omega_1,\\
u^\e = g & \text{on }\Gamma_1,
\end{dcases}
\end{equation}
then one has
\begin{equation}\label{eq:apprx-bdry-C1a}
\sup_{x\in \Omega_\mu} \left| u^\e (x) - \frac{\partial u^\e}{\partial \nu} (0) \inn{x}{\nu} \right| \leq \mu^{1+\alpha},
\end{equation}
where $\nu$ is the inward unit normal to $\Gamma$ at the origin.
\end{lemma}

\begin{proof} 
The idea is similar to Lemma \ref{lemma:apprx-C1a}, whereas the proof here is only more simpler. Let $\mu$ be chosen later. Assume towards a contradiction that for each $k=1,2,\cdots$, there exists a constant $\e_k > 0$ with $\e_k\to 0$, a domain $\Omega_k$ such that $\Gamma_k \in C^2$ contains the origin and the maximal curvature of $\Gamma_{1,k} =\Gamma_k\cap B_1$ is bounded by $1$, a periodic functionals $F_k\in S_0(\lambda,\Lambda)$ on $\cS^n\times\R^n$, and functions $f_k \in  L^\infty(\Omega_{k,1})$ (here and thereafter, we shall denote $\Omega_{k,r} = \Omega_k\cap B_r$ and $\Gamma_{k,r} = \Gamma_k\cap B_r$), $g_k \in C^{1,\alpha}(\Gamma_{k,1} $, $u_k \in C(\overline{\Omega_{k,1}})$ satisfying
\begin{equation}\label{eq:fk-Linf-bdry-C1a}
\norm{f_k}_{L^\infty(\Omega_{k,1})} \leq \e_k,
\end{equation} 
\begin{equation}\label{eq:gk-C1a-bdry-C1a}
g_k(0) = |D_T g_k(0)| = 0,\quad \norm{g_k}_{C^{1,\alpha}(\Gamma_{k,1})} \leq \e_k,
\end{equation} 
\begin{equation}\label{eq:uk-pde-bdry-C1a-apprx}
\begin{dcases}
F_k \left( D^2 u_k,\frac{x}{\e_k} \right) = f_k,\quad |u_k|\leq 1& \text{in }\Omega_{k,1},\\
u_k = g_k & \text{on }\Gamma_{k,1},
\end{dcases}
\end{equation}
such that 
\begin{equation}\label{eq:apprx-bdry-C1a-false}
\sup_{x\in \Omega_{k,\mu}} \left| u_k (x) - \frac{\partial u_k}{\partial \nu_k} (0) \inn{x}{\nu_k} \right| > \mu^{1+\alpha},
\end{equation}
where $\nu_k$ is the inward unit normal to $\Gamma_k$ at the origin. 

Since $\Omega_k\in D(1,\bar\sigma)$, we can assume, after extracting a subsequence if necessary, that $\Omega_k\to \Omega$ for some $\Omega\in D(1,\bar\sigma)$ in the sense of the Hausdorff distance. In particular, $\nu_k \to \nu$ for some $\nu\in \partial B_1$, and $\nu$ is the unit inward normal to $\partial\Omega$ at the origin. Also denoting by $\Phi_k$ and $\Phi$ the rotation mapping associated with $\Omega_k$ and respectively $\Omega$ (as in Definition \ref{definition:domain}) such that $\Phi_k(\nu_k) = e_n = \Phi(\nu)$, we have $\Phi_k \to \Phi$ in $\R^n$.  

Due to the standard boundary $C^{1,\beta}$ estimate \cite[Theorem 1.1]{SS}, there exists $0<\beta<\alpha$ depending only on $n$, $\lambda$ and $\Lambda$ such that $\frac{\partial u_k}{\partial \nu}  \in C^\beta (\Gamma_{k,1/2})$ with
\begin{equation}\label{eq:uk-nu-bdry-C1a-false}
\norm{\frac{\partial u_k}{\partial \nu}}_{C^\beta(\Gamma_{k,1/2})} \leq C_1 \left(\norm{u_k}_{L^\infty(\Omega_{k,1})} + \norm{f_k}_{L^\infty(\Omega_{k,1})} + \norm{g_k}_{C^{1,\alpha}(\Gamma_{k,1})} \right) \leq 2C_1,
\end{equation}
where the second inequality follows from \eqref{eq:fk-Linf-bdry-C1a}, \eqref{eq:gk-C1a-bdry-C1a} and \eqref{eq:uk-pde-bdry-C1a-apprx}. In particular, one has 
\begin{equation}\label{eq:uk-bdry-C1a-false}
\left| u_k(x) - \frac{\partial u_k}{\partial \nu_k}(0) \inn{x}{\nu_k} \right| \leq C_2|x|^{1+\beta},
\end{equation}  
for any $x\in \Omega_{k,1/2}$. Here both $C_1$ and $C_2$ depend only on $n$, $\lambda$, $\Lambda$ and $\alpha$ (note that the maximal curvature of $\Gamma_{k,1}$ is bounded uniformly by $1$). On the other hand, it follows from the standard global H\"older estimate, c.f. \cite[Theorem 2]{S}, that $ u_k \in C_{loc}^\beta (\Gamma_{1,k}\cup\Omega_{1,k})$ and $\| u_k \|_{C^\beta(E)} \leq c_E$ for each $E\Subset \Gamma_{1,k}\cup \Omega_{1,k}$, where $c_E$ depends only on $n$, $\lambda$, $\Lambda$ and $\dist(E,\partial\Omega_{1,k}\setminus \Gamma_{1,k})$, with a possibly smaller $\beta$. Extracting subsequences if necessary, one may assume without loss of generality that \begin{equation}\label{eq:uk-ub-bdry-C1a-false}
u_k\circ\Phi_k^{-1}\circ \Phi \to \bar{u}\quad\text{locally uniformly on }\Gamma_1\cup \Omega_1,
\end{equation} 
for some $\bar{u} \in C_{loc}^\beta (\Gamma_1\cap \Omega_1)$, and 
\begin{equation}\label{eq:uk-nu-ab-bdry-C1a-false}
\frac{\partial u_k}{\partial \nu}(0)  \to \bar{a},
\end{equation}
for some $\bar{a} \in \R$. Hence, taking $k\to \infty$ in \eqref{eq:uk-bdry-C1a-false} we know that $\bar{u}$ satisfies
\begin{equation*}
| \bar{u}(x) - \bar{a} \inn{x}{\nu}| \leq C_2 |x|^{1+\beta},
\end{equation*}
for any $x\in \Omega_{1/2}$, which shows that $\frac{\partial\bar{u}}{\partial\nu}(0)$ exists and 
\begin{equation}\label{eq:ab-ub-nu}
\bar{a} = \frac{\partial \bar{u}}{\partial \nu}(0). 
\end{equation}
Moreover, it follows from \eqref{eq:gk-C1a-bdry-C1a} and \eqref{eq:uk-ub-bdry-C1a-false} that
\begin{equation}\label{eq:ub-bdry-C1a-1}
\norm{\bar{u}}_{L^\infty(\Omega_1)} \leq 1\quad\text{and}\quad\bar{u} = 0 \quad\text{in }\Gamma_1 .
\end{equation}

In addition, arguing as in the proof of Lemma \ref{lemma:apprx-C1a}, one may deduce from the stability result in Lemma \ref{lemma:hom-stab} (here we should apply this lemma for each fixed subdomain of $\Omega_1\cap \Omega_{k,1}$ and then let  $k\to\infty$, so that the stability result holds for any subdomain of $\Omega_1$, and thus itself), the uniform convergence \eqref{eq:uk-ub-bdry-C1a-false} of $u_k\circ\Phi_k^{-1}\circ\Phi\to\bar{u}$ and the assumptions, \eqref{eq:fk-Linf-bdry-C1a} and \eqref{eq:uk-pde-bdry-C1a-apprx}, on $f_k$ and $u_k$ that $\bar{F}_k \to \bar{F}$ locally uniformly on $\cS^n$ for some $\bar{F}\in S_0(\lambda,\Lambda)$, and
\begin{equation}\label{eq:ub-pde-bdry-C1a-apprx}
\bar{F} (D^2 \bar{u}) = 0\quad \text{in }\Omega_1.
\end{equation}
This combined with \eqref{eq:ub-bdry-C1a-1} yields, owing to the standard boundary regularity \cite[Theorem 1.2, Lemma 4.1]{SS}, that
\begin{equation}\label{eq:ub-bdry-C11-apprx}
\left| \bar{u}(x) - \frac{\partial \bar{u}}{\partial \nu}(0) \inn{x}{\nu} \right| \leq \bar{c} |x|^2,
\end{equation} 
for any $x\in \Omega_{1/2}$, where $\bar{c}$ depends only on $n$, $\lambda$, $\Lambda$ and the maximal curvature of $\Gamma_1$. 

Now let us take $0<\mu\leq\frac{1}{2}$ such that $2\bar{c}\mu^{1-\alpha} \leq 1$, which is possible since $0<\alpha<1$. Then we have from \eqref{eq:ub-bdry-C11-apprx} that
\begin{equation}\label{eq:ub-bdry-C11-apprx-re}
\sup_{x\in\Omega_\mu} \left| \bar{u}(x) - \frac{\partial \bar{u}}{\partial \nu}(0) \inn{x}{\nu} \right| \leq \frac{1}{2} \mu^{1+\alpha}.
\end{equation} 
However, this leads us to a contradiction from \eqref{eq:apprx-bdry-C1a-false}, due to \eqref{eq:uk-ub-bdry-C1a-false}, \eqref{eq:uk-nu-ab-bdry-C1a-false} and \eqref{eq:ab-ub-nu}. This finishes the proof. 
\end{proof}

Again it follows an iteration lemma.

\begin{lemma}\label{lemma:iter-bdry-C1a}
Let $\Omega$, $\alpha$, $\mu$ and $\bar\e$ be as in Lemma \ref{lemma:apprx-bdry-C1a}. Then there is $0<\eta\leq 1$, depending only on $n$, $\lambda$, $\Lambda$ and $\alpha$ such that if $\e\leq \bar\e \mu^{k-1}$ for some integer $k\geq 1$, $F\in S_0(\lambda,\Lambda)$ is a periodic functional on $\cS^n\times\R^n$, $f\in L^\infty(\Omega_1)$, $g\in C^{1,\alpha}(\Gamma_1)$ satisfies$g(0) = |D_T g(0)| = 0$, $u^\e$ is a viscosity solution of 
\begin{equation}\label{eq:ue-pde-bdry-C1a}
\begin{dcases}
F \left( D^2 u^\e ,\frac{x}{\e}\right) =f & \text{in }\Omega_1,\\
u^\e = g & \text{on }\Gamma_1,
\end{dcases}
\end{equation}
and the maximal curvature of $\Gamma_1$ is bounded by $\eta$, then there exists $a_k^\e \in \R$ such that
\begin{equation}\label{eq:ake-bdry-C1a}
|a_k^\e|  \leq \frac{n}{\mu(1-\mu^\alpha)} J^\e,\quad\text{and}
\end{equation}
\begin{equation}\label{eq:iter-bdry-C1a}
\sup_{x\in \Omega_{\mu^k}} | u_k^\e (x) - a_k^\e \inn{x}{\nu} | \leq J^\e \mu^{k(1+\alpha)},
\end{equation}
where
\begin{equation}\label{eq:Je-bdry-C1a}
J^\e = \norm{u^\e}_{L^\infty(\Omega_1)} + \frac{1}{\bar\e} \norm{f}_{L^\infty(\Omega_1)} + \frac{2}{\bar\e} \norm{g}_{C^{1,\alpha}(\Gamma_1)}.
\end{equation}
\end{lemma}

\begin{proof}
Now that we have Lemma \ref{lemma:apprx-bdry-C1a}, one may iterate it under the appropriate scaling. Since the argument is very similar to Lemma \ref{lemma:iter-C1a}, we shall leave the details to the reader. 
\end{proof}

We are ready to prove the discrete and uniform boundary estimates. 

\begin{proof}[Proof of Theorem \ref{theorem:bdry-C1a}]
The proof is similar to the proof of Theorem \ref{theorem:int-C1a}. Here we use Lemma \ref{lemma:apprx-bdry-C1a} and Lemma \ref{lemma:iter-bdry-C1a} instead of Lemma \ref{lemma:apprx-C1a} and respectively Lemma \ref{lemma:iter-C1a}. Also note that the smallness condition on the maximal curvature of $\Gamma_1$ in Lemma \ref{lemma:iter-bdry-C1a} can always be achieved by a standard scaling argument. Thus we omit the details to avoid any repeating argument.
\end{proof}


\section{Interior $C^{1,1}$ Estimate}\label{section:int-C11}

This section is devoted to the uniform interior $C^{1,1}$ estimates.

\begin{theorem}\label{theorem:int-C11} 
Let $F\in S_0(\lambda,\Lambda)$ be a periodic functional on $\cS^n\times\R^n$ such that $\bar{F}\in S_2(\lambda,\Lambda,\bar{c},\bar\alpha)$ on $\cS^n$, and $f\in C^\alpha(B_1)$ for some $\alpha\in(0,\bar\alpha)$. Suppose that $u^\e$ is a viscosity solution of 
\begin{equation}\label{eq:main-int-C11}
F\left( D^2 u^\e,\frac{x}{\e}\right) = f \quad\text{in }B_1.
\end{equation}
\begin{enumerate}[(i)]
\item If $F\in S_2(\lambda,\Lambda,\cdot,\cdot)\cap R_0(\cdot,\cdot)$, then for each $x_0\in B_{1/2}$ and any $\e\in(0,\frac{1}{8})$, there exists an affine function $l_{x_0}^\e$, with $l_{x_0}^\e (x_0) = u^\e(x_0)$, such that
\begin{equation}\label{eq:int-C11-e}
|\nabla l_{x_0}^\e| + \sup_{x\in B_{3/4}\setminus \overline{B_\e(x_0)}}\frac{|u^\e (x) - l_{x_0}^\e(x)| }{|x-x_0|^2} \leq C_1 \left( \norm{u^\e}_{L^\infty(B_1)} + \norm{f}_{C^\alpha(B_1)} \right), 
\end{equation} 
where $C_1>0$ depends only on $n$, $\lambda$, $\Lambda$, $\bar{c}$, $\bar\alpha$ and $\alpha$. 
\item Moreover, if $F\in S_2(\lambda,\Lambda,\bar\kappa,\bar\gamma)\cap R_0(\kappa,\gamma)$, then $u^\e \in C^{1,1}(B_{1/2})$ and 
\begin{equation}\label{eq:int-C11}
\norm{u^\e}_{C^{1,1}(B_{1/2})} \leq C_2 \left( \norm{u^\e}_{L^\infty(B_1)} + \norm{f}_{C^\alpha(B_1)} \right),
\end{equation}
where $C_2>0$ depends further on $\bar\kappa$, $\bar\gamma$, $\kappa$ and $\gamma$. 
\end{enumerate}
\end{theorem}

As in Theorem \ref{theorem:int-C1a} (i), the first assertion above amounts to a uniform $C^{1,1}$ estimate for $u^\e$ up to deleted neighborhoods of size $\e$. Again this assertion holds without any control on the modulus of continuity of $F$ in the space variable. This is also new in the context of linear equations (see Remark \ref{remark:iter}). Now the second assertion shows that if $F$ admits $C^{2,\gamma}$ regular solutions in the microscopic scale, then we can fill in the regularity of $u^\e$ in the neighborhoods of size $\e$, completing the full $C^{1,1}$ estimate. 

Let us begin with an approximation lemma. 

\begin{lemma}\label{lemma:apprx} 
Let $\bar\alpha\in(0,1]$ and $\alpha\in(0,\bar\alpha)$ be arbitrary. There are $\mu\in(0,\frac{1}{2}]$, depending only on $\bar{c}$, $\bar\alpha$ and $\alpha$, and $\bar\e\in (0,\frac{1}{2}]$, depending only on $n$, $\lambda$, $\Lambda$, $\bar{c}$, $\bar\alpha$ and $\alpha$, such that for any $\e\leq \bar\e$, any continuous periodic functional $F\in S_0(\lambda,\Lambda)$ on $\cS^n\times\R^n$ having $\bar{F}\in S_2(\lambda,\Lambda,\bar{c},\bar\alpha)$ on $\cS^n$, any function $f\in C^\alpha(B_1)$ and $u^\e\in C(B_1)$ satisfying  
\begin{equation}\label{eq:f-Ca0}
\norm{f}_{C^\alpha(B_1)} \leq \bar\e ,
\end{equation}
\begin{equation}\label{eq:ue-pde}
F\left( D^2 u^\e , \frac{x}{\e}\right) = f \quad \text{in }B_1,
\end{equation}  
\begin{equation}\label{eq:ue-Linf}
u^\e (0) = 0 \quad\text{and}\quad \norm{u^\e}_{L^\infty(B_1)} \leq 1,
\end{equation}
there are some $a^\e \in \R^n$ and $M^\e\in \cS^n$ such that
\begin{equation}\label{eq:Pe}
|a^\e|\leq C_0, \quad |M^\e| \leq \bar{C},\quad \bar{F}( M^\e ) = f(0),\quad\text{and}
\end{equation}
\begin{equation}\label{eq:apprx}
\sup_{x\in B_\mu} \left| u^\e (x) - \inn{a^\e}{x} - \frac{1}{2} \inn{x}{M^\e x} - \e^2 w_F\left( M^\e, \frac{x}{\e}\right) \right| \leq \mu^{2+\alpha},
\end{equation}
where $C_0>0$ depends only on $n$, $\lambda$ and $\Lambda$, and $\bar{C}>0$ depends further on $\bar{c}$.
\end{lemma}

\begin{proof} 
Set $\bar{C} = 2\bar{c} > 0$, and let $C_0$ and $\mu$ be determined later. Suppose towards a contradiction that there is no such $\bar\e$ that the conclusion of this lemma is true. Then there exists a sequence $\{\e_k\}_{k=1}^\infty$ of positive real numbers decreasing to zero, a sequence $\{F_k\}_{k=1}^\infty \subset S_0(\lambda,\Lambda)$ of periodic functionals on $\cS^n\times\R^n$ having $\bar{F}_k \in S_2(\lambda,\Lambda,\bar{c},\bar\alpha)$ on $\cS^n$, and sequences $\{f_k\}_{k=1}^\infty \subset C^\alpha(B_1)$, $\{u_k\}_{k=1}^\infty\subset C(\overline{B_1})$, satisfying  
\begin{equation}\label{eq:fk-Ca0}
\norm{f_k}_{C^\alpha(B_1)}\leq \e_k,
\end{equation}
\begin{equation}\label{eq:uk-pde}
F_k \left( D^2 u_k,\frac{x}{\e_k} \right) = f_k\quad \text{in }B_1,
\end{equation}
\begin{equation}\label{eq:uk-Linf}
u_k(0) = 0 \quad\text{and}\quad \norm{u_k}_{L^\infty(B_1)} \leq 1,
\end{equation}
such that for any $a\in\R^n$ and $M\in\cS^n$ satisfying
\begin{equation}\label{eq:P-false}
|a|\leq \bar{C}, \quad |M|\leq \bar{C}\quad\text{and}\quad \bar{F}_k(M) = f_k(0),
\end{equation}
one has
\begin{equation}\label{eq:apprx-false}
\sup_{x\in B_\mu} \left| u_k (x) - \inn{a}{x} - \frac{1}{2} \inn{x}{M x} - \e_k^2 w_{F_k}\left( M, \frac{x}{\e_k}\right) \right| > \mu^{2+\alpha}.
\end{equation}

Observe from the Krylov-Safanov theory \cite[Proposition 4.9]{CC} and the uniform bound \eqref{eq:uk-Linf} of $u_k$ that $\{u_k\}_{k=1}^\infty$ is uniformly bounded in $C_{loc}^\beta(B_1)$, for some $\beta\in(0,1)$, depending only on $n$, $\lambda$ and $\Lambda$. Hence, one may assume without loss of generality that 
\begin{equation}\label{eq:ubk-11}
u_k\to \bar{u}\quad\text{locally uniformly in }B_1,
\end{equation}
as $k\to\infty$, for some $\bar{u}\in C_{loc}^\beta(B_1)$. Since $F_k \in S_0(\lambda,\Lambda)$ is a periodic functional on $\cS^n\times\R^n$, and $f_k\to 0$ uniformly in $B_1$ as $k\to\infty$, it follows from Lemma \ref{lemma:hom-stab},  \eqref{eq:uk-pde}, \eqref{eq:uk-Linf} and \eqref{eq:ubk-11} that  
\begin{equation}\label{eq:Fbk-Fb}
\bar{F}_k\to \bar{F}\quad\text{locally uniformly on $\cS^n$},
\end{equation}
for some $\bar{F}\in S_0(\lambda,\Lambda)$ on $\cS^n$, along a subsequence, and $\bar{u}$ is a viscosity solution to 
\begin{equation}\label{eq:ubk-111}
\bar{F}(D^2 \bar{u}) = 0\quad\text{in }B_1,\text{ with}
\end{equation}
\begin{equation}\label{eq:ubk-112}
\bar{u}(0) = 0\quad\text{and}\quad \norm{\bar{u}}_{L^\infty(B_1)} \leq 1.
\end{equation}
Moreover, one can assume without losing generality that \eqref{eq:Fbk-Fb} holds in the full sequence. 

On the other hand, from $\bar{F}_k\in S_0(\lambda,\Lambda)$ along with the Krylov theory \cite[Corollary 5.7]{CC} we also know that $\bar{F}_k\in S_1(\lambda,\Lambda,c_0,\alpha_0)$, for some $c_0>0$ and $\alpha_0\in(0,1)$, both depending only on $n$, $\lambda$ and $\Lambda$. Thus, we can apply the almost $C^{1,\frac{\alpha_0}{2}}$ estimate (Theorem \ref{theorem:int-C1a} (i)) to $u_k$, from which we obtain a vector $a_k \in \R^n$ such that 
\begin{equation}\label{eq:uk-int-C1a-almst}
|a_k| + \sup_{x\in B_{1/2}\setminus \overline{B_{\e_k}}} \frac{|u_k (x) - \inn{a_k}{x}|}{|x|^{1+\frac{\alpha_0}{2}}} \leq C_0,
\end{equation} 
where $C_0>0$ depends only on $n$, $\lambda$, $\Lambda$, $c_0$ and $\alpha_0$, hence on the first three parameters only. 

As $\{a_k\}_{k=1}^\infty$ being a bounded sequence in $\R^n$, one can extract a subsequence $\{a_{k_i}\}_{i=1}^\infty$ of $\{a_k\}_{k=1}^\infty$ such that $a_{k_i} \to \bar{a}$, as $i\to\infty$, for certain $\bar{a}\in\R^n$. Passing to the limit in \eqref{eq:uk-int-C1a-almst}, we observe from the uniform convergence \eqref{eq:ubk-11} of $u_{k_i}\to \bar{u}$ over $B_{1/2}$ that
\begin{equation}\label{eq:ub-int-C1a}
|\bar{a}| + \sup_{x\in B_{1/2} \setminus\{0\}} \frac{ | \bar{u}(x) - \inn{\bar{a}}{x} | }{|x|^{1+\frac{\alpha_0}{2}}} \leq C_0. 
\end{equation} 
However, since $\bar{F}\in S_0(\lambda,\Lambda)$, we also have $\bar{F}\in S_1(\lambda,\Lambda,c_0,\alpha_0)$, which along with \eqref{eq:ub-int-C1a} implies that $\bar{a} = \nabla \bar{u}(0)$. Therefore, we deduce that 
\begin{equation}\label{eq:ub-int-C1a-2}
a_k \to \bar{a},
\end{equation}
as $k\to\infty$ along the full sequence. 

Now since we assume $\bar{F}_k \in S_2(\lambda,\Lambda,\bar{c},\bar\alpha)$ for each $k\in\N$, Lemma \ref{lemma:limit} together with \eqref{eq:Fbk-Fb} yields that $\bar{F} \in S_2(\lambda,\Lambda,\bar{c},\bar\alpha)$. Noting that $\bar{u}$ is a viscosity solution of \eqref{eq:ubk-111}, it follows from the definition of the class $S_2$ that 
\begin{equation}\label{eq:ub-C2ab}
\norm{\bar{u}}_{C^{2,\alpha}(\bar{B}_{1/2})} \leq \bar{c}.
\end{equation}
In particular, with $\bar{M} = D^2 \bar{u}(0)$ (and $\bar{a} = \nabla \bar{u}(0)$), we have 
\begin{equation}\label{eq:Mb}
|\bar{M}| \leq \bar{c}\quad\text{and}\quad \sup_{x\in B_{1/2} \setminus\{0\}} \frac{| \bar{u}(x) - \inn{\bar{a}}{ x} - \frac{1}{2} \inn{x}{\bar{M} x} |}{|x|^{2+\alpha}} \leq \bar{c}. 
\end{equation} 
Thus, selecting a sufficiently small $\mu\leq\frac{1}{2}$ such that 
\begin{equation}\label{eq:mu-int}
\bar{c}\mu^{\bar\alpha} \leq\frac{1}{2} \mu^{\alpha},
\end{equation}
we arrive at 
\begin{equation}\label{eq:ub-apprx}
\sup_{x\in B_\mu} \left| \bar{u}(x) - \inn{\bar{a}}{x} - \frac{1}{2} \inn{x}{\bar{M} x} \right| \leq \frac{1}{2} \mu^{2+\alpha}. 
\end{equation}
Clearly, the smallness of $\mu$ depends only on $\bar\alpha$, $\alpha$ and $\bar{c}$. Let us also remark that as $\bar{u}$ being a classical solution of \eqref{eq:ubk-111} in $B_{1/2}$, one has  
\begin{equation}\label{eq:ub-apprx3}
\bar{F}(\bar{M}) = 0.
\end{equation}

We are going to find a sequence $\{M_k\}_{k=1}^\infty\subset\cS^n$ converging to $M$ and satisfying \eqref{eq:P-false}.
Setting
\begin{equation*}
\delta_k = |\bar{F}_k(\bar{M})|
\end{equation*}
we know from \eqref{eq:Fbk-Fb} and\eqref{eq:ub-apprx3} $\delta_k\to 0$, so it follows from the property $\bar{F}_k\in S_0(\lambda,\Lambda)$ and the assumption \eqref{eq:fk-Ca0} that 
\begin{equation*}
\bar{F}_k(\bar{M} + \lambda^{-1}(\e_k + \delta_k) I) \geq \bar{F}_k(\bar{M}) + \delta_k + \e_k  \geq \e_k \geq |f_k(0)|,
\end{equation*} 
where in the last inequality we used \eqref{eq:ub-apprx3}. Similarly, one has
\begin{equation*}
\bar{F}_k(\bar{M} - \Lambda^{-1}(\e_k + \delta_k)I) \leq \bar{F}_k(\bar{M}) -\delta_k - \e_k \leq f_k(0). 
\end{equation*}
Thus, by the intermediate value theorem, there must exists $M_k\in\cS^n$, for each $k\in\N$, satisfying 
\begin{equation}\label{eq:Mk-fk}
\bar{F}_k(M_k) = f_k(0)\quad\text{and}\quad |M_k - \bar{M}| \leq \lambda^{-1} (\e_k + \delta_k).
\end{equation}
Especially, $M_k \to \bar{M}$ as $k\to\infty$, and it follows from \eqref{eq:Mb} that
\begin{equation}\label{eq:Mk-Linf}
|M_k| \leq \bar{C},
\end{equation}
for all sufficiently large $k$'s, since we have taken $\bar{C} = 2\bar{c}$ from the beginning. Therefore, for all large $k$'s, $a_k$ and $M_k$ satisfies \eqref{eq:P-false}, so we deduce from \eqref{eq:apprx-false} that
\begin{equation}\label{eq:Mk-false}
\sup_{x\in B_\mu} \left| u_k (x)  - \inn{a_k}{x} - \frac{1}{2} \inn{x}{M_k x} - \e_k^2 w_{F_k}\left( M_k, \frac{x}{\e_k}\right) \right| > \mu^{2+\alpha}.
\end{equation}

Due to \eqref{eq:wF-Ca} and \eqref{eq:Mk-Linf}, 
\begin{equation}\label{eq:wFk-Mk}
\norm{w_{F_k}(M_k,\cdot)}_{L^\infty(\R^n)} \leq C_1|M_k|\leq 2C_1\bar{c},
\end{equation}
where $C_1$ depends only on $n$, $\lambda$ and $\Lambda$. Hence, from \eqref{eq:ubk-11}, \eqref{eq:ub-int-C1a-2}, \eqref{eq:Mk-false} and \eqref{eq:Mk-Linf}, passing to the limit in \eqref{eq:Mk-false} yields that 
\begin{equation*}
\sup_{x\in \bar{B}_\mu} \left| \bar{u} (x) - \inn{\bar{a}}{x} - \frac{1}{2} \inn{x}{\bar{M} x} \right| \geq \mu^{2+\alpha},
\end{equation*}
a contradiction to \eqref{eq:ub-apprx}. Thus, there must exist $\bar\e\in (0,\frac{1}{2}]$ such that the assertion of this lemma is true. Moreover, the smallness of $\bar\e $ is determined only by the fixed parameters chosen from the beginning, which are precisely $n$, $\lambda$, $\Lambda$, $\bar{c}$ and $\bar\alpha$. This finishes the proof.
\end{proof}

Next we establish an iteration lemma.

\begin{lemma}\label{lemma:iter} 
Let $\alpha$, $\mu$, $\bar\e$, $C_0$ and $\bar{C}$ be as in Lemma \ref{lemma:apprx}. If $\e \leq \bar\e \mu^{k-1}$ for some integer $k\geq 1$, $F\in S_2(\lambda,\Lambda,\cdot,\cdot)\cap R_0(\cdot,\cdot)$ is a periodic functional on $\cS^n\times\R^n$ having $\bar{F} \in S_2(\lambda,\Lambda,\bar{c},\bar\alpha)$ on $\cS^n$, $f\in C^\alpha(B_1)$ and $u^\e$ is a viscosity solution of 
\begin{equation}\label{eq:ue-pde-gen}
F \left( D^2 u^\e,\frac{x}{\e}\right) = f\quad\text{in }B_1,
\end{equation}
then for each $l\in \{1,\cdots,k\}$, there exist $a_l^\e\in\R^n$ and $M_l^\e\in\cS^n$, such that $a_0^\e = 0$, $M_0^\e = 0$, 
\begin{equation}\label{eq:ake}
|a_l^\e - a_{l-1}^\e| \leq C_0 J^\e  \mu^{(l-1)(1+\alpha)},
\end{equation}
\begin{equation}\label{eq:Mke}
|M_l^\e - M_{l-1}^\e | \leq \bar{C} J^\e \mu^{(l-1)\alpha}\quad\text{with}\quad \bar{F}(M_l^\e) = f(0),\quad\text{and}
\end{equation}
\begin{equation}\label{eq:iter}
\sup_{x\in B_{\mu^l}} \left| u^\e (x) - u^\e(0) - \inn{a_l^\e}{x} - \frac{1}{2} \inn{x}{M_l^\e x} - \e^2 w_F \left( M_l^\e , \frac{x}{\e} \right) \right| \leq J^\e \mu^{l(2+\alpha)},
\end{equation}
where 
\begin{equation}\label{eq:Ke}
J^\e = \norm{u^\e}_{L^\infty(B_1)} + \frac{1}{\bar\e } \norm{f}_{C^\alpha(B_1)}.
\end{equation}
\end{lemma}

\begin{proof} 
As in the proof of Lemma \ref{lemma:iter-C1a}, we shall assume that Lemma \ref{lemma:iter} holds with $\e\leq \bar\e \mu^{k-1}$ for some $k\geq 1$, and prove that it continues to hold when $\e\leq \bar\e\mu^k$. Thus, we have $a_l^\e \in \R^n$ and $M_l^\e \in\cS^n$ for $l\in \{1,\cdots,k\}$, with $|a_0^\e | = |M_0^\e| = 0$, satisfying \eqref{eq:ake} as well as \eqref{eq:Mke}, and construct $a_{k+1}^\e\in\R^n$ and $M_{k+1}^\e\in\cS^n$ that verify \eqref{eq:ake}, \eqref{eq:Mke} and \eqref{eq:iter} with $l$ replaced by $k+1$. 

Define, for $x\in B_1$, 
\begin{equation}\label{eq:vke}
u_k^\e(x)  = \frac{u^\e(\mu^k x)  - u^\e(0) - \mu^k \inn{a_k^\e}{x} - \frac{1}{2}\mu^{2k} x^T M_k^\e x - \e^2 w_F ( M_k^\e, \e^{-1}\mu^k x)}{J^\e \mu^{k(2+\alpha)}},
\end{equation}
and 
\begin{equation}\label{eq:vpk}
f_k^\e (x) = \frac{f (\mu^k x) - f (0)}{J^\e \mu^{k\alpha}}.
\end{equation}
Also set, for $M\in\cS^n$ and $y\in\R^n$, 
\begin{equation}\label{eq:Gke}
F_k^\e (N,y) = \frac{F(J^\e \mu^{k\alpha}  N + M_k^\e + D_y^2 w_F(M_k^\e, y), y) - f(0)}{J^\e \mu^{k\alpha}}.
\end{equation}
which is clearly a uniformly elliptic, continuous and periodic functional on $\cS^n\times\R^n$, with ellipticity constants $\lambda$ and $\Lambda$; here the continuity of $F_k^\e$ follows from $D_y^2 w_F(M_k^\e,\cdot)\in C(\R^n)$, which is again ensured by the assumption $F\in S_2(\lambda,\Lambda,\cdot,\cdot)\cap R_0(\cdot,\cdot)$. As $u^\e$ being a viscosity solution to \eqref{eq:ue-pde-gen} and $w_F(M_k^\e,\cdot)\in C^2(\R^n)$, $u_k^\e$ becomes a viscosity solution of 
\begin{equation}\label{eq:vke-pde-iter}
F_k^\e \left( D^2 u_k^\e,\frac{\mu^k x}{\e}\right) = f_k^\e\quad \text{in }B_1,\text{ with}
\end{equation}
\begin{equation}\label{e:vke-Linf-iter}
\norm{u_k^\e}_{L^\infty(B_1)} \leq 1,
\end{equation}
where the last inequality follows from the induction hypothesis \eqref{eq:iter}. Moreover, it is clear from the choice \eqref{eq:Ke} of $J^\e$ that 
\begin{equation}\label{eq:fke-Ca0}
\norm{f_k^\e}_{C^\alpha(B_1)} \leq \frac{1}{J^\e} \norm{f}_{C^\alpha(B_{\mu^k})} \leq \bar\e . 
\end{equation} 

From the last identity in the induction hypothesis \eqref{eq:ake}, and the definition of the effective functional $\bar{F}$, we know that 
\begin{equation}\label{eq:Fke-re}
\begin{split}
F_k^\e (N,y) &= \frac{F(J^\e \mu^{k\alpha}  N + M_k^\e + D_y^2 w_F(M_k^\e, y), y) - \bar{F}(M_k^\e)}{J^\e \mu^{k\alpha}}\\
& = \frac{F(J^\e \mu^{k\alpha}  N + M_k^\e + D_y^2 w_F(M_k^\e, y), y) - F( M_k^\e + D_y^2 w_F(M_k^\e,y),y)}{J^\e \mu^{k\alpha}},
\end{split}
\end{equation}
so $F_k$ has zero source term. This implies that $F_k^\e \in S_0(\lambda,\Lambda)$ on $\cS^n\times\R^n$.

Let us denote by $\bar{F}_k^\e$ the effective functional corresponding to $F_k^\e$. Following the notation in Lemma \ref{lemma:GF}, we have $F_k^\e = F_{M_k^\e, J^\e \mu^{k\alpha}}$, so 
\begin{equation}\label{eq:Gbke}
\bar{F}_k^\e(N) = \frac{\bar{F} ( \mu^{k\alpha}N + M_k^\e ) - \bar{F} (M_k^\e)}{\mu^{k\alpha}},
\end{equation}
and
\begin{equation}\label{eq:wFke}
w_{F_k^\e}(N,y) = \frac{w(\mu^{k\alpha}N + M_k^\e,y) - w(M_k^\e,y)}{\mu^{k\alpha}},
\end{equation}
for any $(N,y)\in\cS^n\times\R^n$. In particular, we have $\bar{F}_k^\e \in S_2(\lambda,\Lambda,\bar{c},\bar\alpha)$ as well. 

As $F_k^\e \in S_0(\lambda,\Lambda)$ being a periodic functional with $\bar{F}_k^\e \in S_2(\lambda,\Lambda,\bar{c},\bar\alpha)$, $f_k^\e \in C^\alpha(B_1)$ with the smallness condition \eqref{eq:fke-Ca0}, and $\e\mu^{-k} \leq \bar\e$, we can apply the approximation lemma (Lemma \ref{lemma:apprx}) to the normalized equation \eqref{eq:vke-pde-iter}. This yields a vector $b_k^\e \in \R^n$ and a matrix $N_k^\e \in \cS^n$, satisfying
\begin{equation}\label{eq:Nke}
|b_k^\e| \leq C_0,\quad |N_k^\e| \leq \bar{C}\quad\text{and}\quad \bar{F}_k^\e(N_k^\e) = f_k^\e(0) = 0,
\end{equation} 
where $C_0$ and $\bar{C}$ are the same constants appearing in \eqref{eq:Pe}, such that 
\begin{equation}\label{eq:vke-apprx}
\sup_{x\in B_\mu} \left| u_k^\e (x) - \inn{b_k^\e}{x} - \frac{1}{2} \inn{x}{N_k^\e x} - \frac{\e^2}{\mu^{2k}} w_{F_k^\e} \left( N_k^\e , \frac{\mu^k x}{\e} \right) \right| \leq \mu^{2+\alpha},
\end{equation}

Therefore, setting
\begin{equation}\label{eq:ak1e}
a_{k+1}^\e = a_k^\e + J^\e \mu^{k(1+\alpha)} b_k^\e,
\end{equation}
and
\begin{equation}\label{eq:Mk1e}
M_{k+1}^\e = M_k^\e + J^\e \mu^{k\alpha} N_k^\e,
\end{equation}
we can rephrase \eqref{eq:vke-apprx}, through \eqref{eq:Gbke} and \eqref{eq:wFke}. as 
\begin{equation*}
\begin{split}
& \sup_{x\in B_{\mu^{k+1}}} \left| u^\e (x) - u^\e (0) - \inn{a_{k+1}^\e}{x} - \frac{1}{2} \inn{x}{M_{k+1}^\e x} - \e^2 w_F \left( M_{k+1}^\e,\frac{x}{\e}\right)  \right| \\
& \leq J^\e\mu^{(k+1)(2+\alpha)}.
\end{split}
\end{equation*}
Due to \eqref{eq:Nke} and the induction hypothesis \eqref{eq:ake}, we have 
\begin{equation*}
|a_{k+1}^\e - a_k^\e| \leq C_0 J^\e \mu^{k(1+\alpha)},
\end{equation*}
and from \eqref{eq:Mke} we also derive that 
\begin{equation*}
|M_{k+1}^\e - M_k^\e| \leq \bar{C} J^\e \mu^{k\alpha}.
\end{equation*}
Moreover, it follows from \eqref{eq:Gbke} that
\begin{equation*}
\bar{F}(M_{k+1}^\e) = \bar{F}(M_k^\e) + \mu^{k\alpha} \bar{F}_k^\e(N_k^\e) = f(0). 
\end{equation*} 
Therefore, both \eqref{eq:ake} and \eqref{eq:Mke} are verified with $l$ replaced by $k+1$, which finishes the proof owing to the induction principle. 
\end{proof}

\begin{remark}\label{remark:iter}
One can see that if $F$ were a linear functional so that $F(M,y) = \tr(A(y)M)$ for some uniformly elliptic matrix $A(y)$, then $F_k^\e$ defined as in \eqref{eq:Fke-re} becomes
\begin{equation}
F_k^\e (N,y) = \tr(A(y) N),
\end{equation}
whence there is no need to use $C^2$ regular interior corrector to construct the equation for $u_k^\e$. This is the main reason why it is enough to assume uniform ellipticity only to work with the argument above, in the framework of linear equations. 
\end{remark}

Let us prove Theorem \ref{theorem:int-C11}. 

\begin{proof}[Proof of Theorem \ref{theorem:int-C11}]
With the iteration lemma (Lemma \ref{lemma:iter}), we can prove Theorem \ref{theorem:int-C11} by following the argument in the proof of Theorem \ref{theorem:int-C1a} as well as that of \cite[Theorem 1 (iii)]{AL2}. We omit the details. 
\end{proof}


\section{Boundary $C^{1,1}$ Estimate}\label{section:bdry-C11}

This section is devoted to the uniform boundary $C^{1,1}$ estimate. 

\begin{theorem}\label{theorem:bdry-C11}
Let $\Omega\in D(\bar\tau,\bar\sigma)$, and suppose that $F\in S_2(\lambda,\Lambda,\bar\kappa,\bar\gamma)\cap R_1(\kappa,\gamma)$ is a periodic functional on $\cS^n\times\R^n$ having $\bar{F} \in S_2(\lambda,\Lambda,\bar{c},\bar\alpha)$ on $\cS^n$, $f\in C^\alpha(\Omega_1)$, $g\in C^{2,\alpha}(\Gamma_1)$, with some $0<\alpha<\bar\alpha$, and $u^\e$ is a viscosity solution of 
\begin{equation}\label{eq:main-bdry-C11}
\begin{dcases}
F \left( D^2 u^\e,\frac{x}{\e}\right) = f & \text{in }\Omega_1,\\
u^\e = g & \text{on }\Gamma_1.
\end{dcases}
\end{equation}
Then $u^\e \in C^{1,1}(\Omega_{1/2})$ and 
\begin{equation}\label{eq:bdry-C11}
\norm{u^\e}_{C^{1,1}(\Omega_{1/2})} \leq C,
\end{equation}
where $C>0$ depends only on $n$, $\lambda$, $\Lambda$, $\bar\kappa$, $\bar\gamma$, $\kappa$, $\gamma$, $\bar{c}$, $\bar\alpha$, $\alpha$, $\bar\tau$, $\bar\sigma$, $\norm{f}_{C^\alpha(\Omega_1)}$, $\norm{g}_{C^{2,\alpha}(\Gamma_1)}$ and $\sup_{\e>0} \norm{u^\e}_{L^\infty(\Omega_1)}$. 
\end{theorem}

We have already seen in the iteration scheme for the interior $C^{1,1}$ estimate that the nonlinearity of the governing equation results a new equation at each iterative step, which amounts to the effect coming from the correction in the previous steps; see the functional $F_k^\e$ defined in \eqref{eq:Gke} that appears in the $k$-th iteration step. In the case of the boundary estimates, one has to involve a boundary layer corrector, and the same phenomenon occurs. However, the problem becomes more difficult, since the boundary layer corrector does not oscillate in the periodic manner, and since it also solves a rapidly oscillating nonlinear equation. It should be stressed that we do not encounter this issue in the context of linear equations, since the linearity annihilates the effect coming from the previous correction; see Remark \ref{remark:iter} for instance. 

Let us begin with a uniform boundary $W^{2,p}$ estimates, for any $p>n$. This can be understood as a byproduct of the uniform interior $C^{1,1}$ estimates (Theorem \ref{theorem:int-C11}) and the boundary $C^{1,\alpha}$ estimates (Theorem \ref{theorem:bdry-C1a}), for any $\alpha\in(0,1)$.

\begin{proposition}\label{proposition:bdry-W2p}
Let $\Omega$ be a domain with boundary $\Gamma\in C^2$ containing the origin. Suppose that $F \in S_2 (\lambda,\Lambda,\bar\kappa,\bar\gamma) \cap R_0(\kappa,\gamma)$ is a periodic functional on $\cS^n\times\R^n$ having $\bar{F} \in S_1(\lambda,\Lambda,\bar{c},1)$ on $\cS^n$, $f\in C^\alpha(\Omega_1)$ for some $\alpha\in(0,1)$, $g\in C^{1,1}(\Gamma_1)$ and $u^\e$ is a viscosity solution of 
\begin{equation}\label{eq:main-bdry-W2p}
\begin{dcases}
F \left( D^2 u^\e, \frac{x}{\e} \right) = f & \text{in }\Omega_1,\\
u^\e = g& \text{on }\Gamma_1.
\end{dcases}
\end{equation}
Let $p>0$ be any. Then $u^\e \in W^{2,p}(\Omega_{1/2})$ and 
\begin{equation}\label{eq:bdry-W2p}
\norm{u^\e}_{W^{2,p}(\Omega_{1/2})} \leq C \left( \norm{u^\e}_{L^\infty(\Omega_1)} + \norm{f}_{C^\alpha(\Omega_1)} + \norm{g}_{C^{1,1}(\Gamma_1)} \right), 
\end{equation} 
where $C>0$ depends only on $n$, $\lambda$, $\Lambda$, $\bar\kappa$, $\bar\gamma$, $\kappa$, $\gamma$, $\bar{c}$, $\alpha$, the maximal curvature of $\Gamma_1$ and $p$. 
\end{proposition}

\begin{proof}
To simplify the notation, let us call $\delta(x)$ the distance function $d(x,\Gamma_1)$. Fix $x_0\in \Omega_{1/2}$ and consider an auxiliary function $v^\e$ on $B_{\delta(x_0)}(x_0)$ defined by 
\begin{equation*}
v^\e (x) = u^\e (x) - u^\e (x_0) - Du^\e (x_0)\cdot(x-x_0). 
\end{equation*}
In view of \eqref{eq:main-bdry-W2p}, $v^\e$ solves 
\begin{equation*}
F \left( D^2 v^\e,\frac{x}{\e}\right) = f \quad\text{in }\Omega_1,
\end{equation*}
and since $F$ and $f$ satisfy the assumption of the uniform interior $C^{1,1}$ estimate (Theorem \ref{theorem:int-C11}), it follows from a scaled version of \eqref{eq:int-C11} that 
\begin{equation}\label{eq:ve-C11-bdry-W2p}
|D^2 u^\e (x_0)| = |D^2 v^\e(x_0)| \leq C_1 \left( \frac{\norm{v^\e}_{L^\infty(B_{\delta(x_0)}(x_0))}}{\delta(x_0)^2} + \norm{f}_{C^\alpha(B_{\delta(x_0)}(x_0))} \right),
\end{equation}
where $C_1$ depends only on $n$, $\lambda$, $\Lambda$, $\bar\kappa$, $\bar\gamma$, $\kappa$, $\gamma$ and $\alpha$. 

In order to estimate the $L^\infty$ norm of $v^\e$ on $B_{\delta(x_0)}(x_0)$, we consider a uniform boundary $C^{1,\alpha}$ estimate for $u^\e$. From the assumption that $\bar{F} \in S_1(\lambda,\Lambda,\bar{c},1)$, one can apply Theorem \ref{theorem:int-C1a} (ii) to derive that $u^\e \in C^{1,\alpha}(\Omega_{3/4})$, for any $0<\alpha<1$, and 
\begin{equation}\label{eq:ue-bdry-C1a-bdry-W2p}
\norm{u^\e}_{C^{1,\alpha}(\Omega_{3/4})} \leq C_2 \left( \norm{u^\e}_{L^\infty(\Omega_1)} + \norm{f}_{C^\alpha(\Omega_1)} + \norm{g}_{C^{1,1}(\Gamma_1)} \right),
\end{equation}
where $C_2$ depends only on  $n$, $\lambda$, $\Lambda$, $\bar\kappa$, $\bar\gamma$, $\kappa$, $\gamma$, $\bar{c}$ and $\alpha$. Especially this implies that 
\begin{equation*}
\norm{v^\e}_{L^\infty(B_{\delta(x_0)}(x_0))} \leq C_2 \delta(x_0)^{1+\alpha}\left( \norm{u^\e}_{L^\infty(\Omega_1)} + \norm{f}_{C^{0,1}(\Omega_1)} +  \norm{g}_{C^{1,1}(\Gamma_1)} \right),
\end{equation*}
which in turn yields in \eqref{eq:ve-C11-bdry-W2p} that
\begin{equation}\label{eq:ue-C11-bdry-W2p}
|D^2 u^\e (x_0)| \leq C_3 \delta(x_0)^{\alpha-1} \left( \norm{u^\e}_{L^\infty(\Omega_1)} +\norm{f}_{C^\alpha(\Omega_1)} + \norm{g}_{C^{1,1}(\Gamma_1)}\right),  
\end{equation}
where $C_3$ depends on the same parameters that determine both $C_1$ and $C_2$. 

Now given $p>0$, let us choose $\alpha$ as close as $1$ from below such that $(1-\alpha)p<1$. Then we have 
\begin{equation*}
\int_{\Omega_{1/2}} \delta(x)^{(\alpha-1)p} \,dx \leq C_4,
\end{equation*}
where $C_4$ depends only on the maximal curvature of $\Gamma_1$ and $p$. This together with \eqref{eq:ue-C11-bdry-W2p} implies that 
\begin{equation}
\norm{D^2 u^\e}_{L^p(\Omega_{1/2})} \leq C_5 \left( \norm{u^\e}_{L^\infty(\Omega_1)} + \norm{f}_{C^\alpha(\Omega_1)} +  \norm{g}_{C^{1,1}(\Gamma_1)}\right), 
\end{equation} 
where $C_5$ depends only on $n$, $\lambda$, $\Lambda$, $\bar\kappa$, $\bar\gamma$, $\kappa$, $\gamma$, $\bar{c}$, the maximal curvature of $\Gamma_1$ and $p$. Hence, the proof is finished.
\end{proof}

As in the previous sections, let us present an approximation lemma. 

\begin{lemma}\label{lemma:apprx-bdry} 
Let $\Omega \in D(1,\bar\sigma)$. Then there are $0<\bar\alpha<\bar\sigma$ and $\bar{C}>0$, depending only on $n$, $\lambda$, $\Lambda$ and $\bar\sigma$, such that the following is true: for any $0<\alpha<\bar\alpha$, $L>0$ and $p>n$, one can choose $0<\mu\leq\frac{1}{2}$, depending only on depending only on $n$, $\lambda$, $\Lambda$, $\bar\sigma$ and $\alpha$, and $0<\bar\e\leq\frac{1}{2}$, depending only on $n$, $\lambda$, $\Lambda$, $\bar\kappa$, $\bar\gamma$, $\kappa$, $\gamma$, $\bar\sigma$, $\alpha$, $L$ and $p$, such that for any $\e\leq \bar\e$, any periodic functional $F\in S_2(\lambda,\Lambda,\bar\kappa,\bar\gamma)\cap R_0(\kappa,\gamma)$ on $\cS^n\times\R^n$, any periodic functional $G\in S_0(\lambda,\Lambda)$ on $\cS^n\times\Omega\times\R^n$, satisfying for any $(N,y)\in\cS^n\times\R^n$,
\begin{equation}\label{eq:G-F-Lp-bdry}
\norm{G(N,\cdot,y) - F(N,y)}_{L^p(\Omega_1)} \leq |N| L,\quad\text{and}
\end{equation}
\begin{equation}\label{eq:G-F-Linf-bdry}
\norm{ d(\cdot,\Gamma_1)^2 (G(N,\cdot,y) - F(N,y))}_{L^\infty(\Omega_1)} \leq |N|\bar\e,
\end{equation} 
and any $f\in C^\alpha(\Omega_1)$, $g\in C^{2,\alpha}(\Gamma_1)$ and $u^\e \in C(\overline{\Omega_1})$, satisfying
\begin{equation}\label{eq:f-Ca0-bdry}
\norm{f}_{C^\alpha(\Omega_1)} \leq \bar\e,
\end{equation}
\begin{equation}\label{eq:g-C2a0-bdry}
g(0) = |D_T g(0)| = 0,\quad \norm{g}_{C^{2,\alpha}(\Gamma_1)} \leq \bar\e,\quad\text{and}
\end{equation}
\begin{equation}\label{eq:ue-pde-bdry}
\begin{dcases}
G \left( D^2 u^\e,x,\frac{x}{\e}\right) = f,\quad |u^\e|\leq 1 & \text{in }\Omega\cap B_1,\\
u^\e = g & \text{on }\partial\Omega\cap B_1,
\end{dcases}
\end{equation}
there exist a matrix $M^\e\in\cS^n$ such that 
\begin{equation}\label{eq:Me-bdry}
|M^\e|\leq \bar{C},\quad (I - \nu\otimes\nu) M^\e = D_T^2 g(0),\quad \bar{F}(M^\e) = f(0),\quad\text{and}
\end{equation}
\begin{equation}\label{eq:apprx-bdry}
\sup_{x\in \Omega_\mu} \left| u^\e (x) - \frac{\partial u^\e}{\partial\nu} (0)\inn{x}{\nu} - \frac{1}{2} \inn{x}{M^\e x} - \e^2 w_F \left( M^\e,\frac{x}{\e}\right) - v^\e (x) \right|\leq \mu^{2+\alpha},
\end{equation}
where $v^\e$ is the viscosity solution of 
\begin{equation}\label{eq:ve-pde-bdry}
\begin{dcases}
G \left( M^\e + D_y^2 w_F\left(M^\e,\frac{x}{\e}\right) + D_x^2 v^\e, x, \frac{x}{\e}\right) = \bar{F}(M^\e) = f(0)& \text{in }\Omega_1,\\
v^\e = - \e^2 w_F \left( M^\e,\frac{x}{\e} \right) & \text{on }\partial\Omega_1. 
\end{dcases}
\end{equation}
\end{lemma}

\begin{proof}
Let $\bar\alpha$ and $\bar{c}$ be determined later, set $\bar{C} = 2\bar{c}$ and fix $0<\alpha<\bar\alpha$. Also set $\mu$ to a small positive constant to be determined later. Assume to the contrary that there is no such $\bar\e $ so that the conclusion of this lemma holds. Then for each $k=1,2,\cdots$, one can choose a constant $\e_k>0$ with $\e_k\to 0$, a bounded domain $\Omega_k \in D(1,\bar\sigma)$ with $\nu_k$ being the inward unit normal at the origin (we shall write $\Omega_{k,r} = \Omega_k\cap B_r$ and $\Gamma_{k,r} = \Gamma_k\cap B_r$ throughout this proof), a periodic functional $F_k \in S_2 (\lambda,\Lambda,\bar\kappa,\bar\gamma)\cap R_0(\kappa,\gamma)$ on $\cS^n\times\R^n$, another functional $G_k\in S_0(\lambda,\Lambda)$ on $\cS^n\times\Omega_{k,1} \times \R^n$ such that $G_k(N,x,\cdot)$ is periodic on $\R^n$ for each $(N,x)\in\cS^n\times\Omega_{k,1}$, and thatfor any $(N,y)\in\cS^n\times\R^n$, 
\begin{equation}\label{eq:Gk-Fk-Lp-bdry}
\norm{G_k(N,\cdot,y) - F_k(N,y)}_{L^p(\Omega_{k,1})} \leq |N| L,\quad\text{and}
\end{equation}
\begin{equation}\label{eq:Gk-Fk-Linf-bdry}
\norm{ d(\cdot,\Gamma_{k,1})^2 (G_k(N,\cdot,y) - F_k(N,y))}_{L^\infty(\Omega_{k,1})} \leq |N|\e_k,
\end{equation}
and functions $f_k \in  C^\alpha(\Omega_{k,1})$, $g_k \in C^{2,\alpha}(\Gamma_{k,1})$, $u_k \in C(\overline{\Omega_{k,1}})$ satisfying 
\begin{equation}\label{eq:fk-Ca0-bdry}
\norm{f_k}_{C^\alpha(\Omega_{k,1})} \leq \e_k,
\end{equation}
\begin{equation}\label{eq:gk-C2a0-bdry}
g_k(0) = |D_T g_k(0)| = 0, \quad \norm{g_k}_{C^{2,\alpha}(\Gamma_{k,1})} \leq \e_k,\quad\text{and}
\end{equation}
\begin{equation}\label{eq:uk-pde-bdry}
\begin{dcases}
G_k\left(D^2 u_k, x, \frac{x}{\e_k} \right) = f_k,\quad |u_k|\leq 1 & \text{in }\Omega_{k,1},\\
u_k = g_k & \text{on }\Gamma_{k,1},
\end{dcases}
\end{equation}
such that for any $M\in\cS^n$ satisfying 
\begin{equation}\label{eq:M-bdry-false}
|M|\leq \bar{C}, \quad (I - \nu_k\otimes \nu_k) M = D_T^2 g_k(0), \quad \bar{F}_k(M) = f_k(0),
\end{equation}
one has 
\begin{equation}\label{eq:apprx-bdry-false}
\sup_{x\in \Omega_{k,\mu}} \left| u_k (x) - \frac{\partial u_k}{\partial\nu_k}(0) \inn{x}{\nu_k} - \frac{1}{2} \inn{x}{M x} - \e_k^2 w_{F_k} \left(M ,\frac{x}{\e_k}\right) - v_k (x)\right| > \mu^{2+\alpha},
\end{equation}
where $v_k\in C(\overline{\Omega_{k,1}})$ is the viscosity solution of 
\begin{equation}\label{eq:vk-bdry-false}
\begin{dcases}
G_k \left( M + D_y^2 w_{F_k}\left(M,\frac{x}{\e_k}\right) + D_x^2 v_k, x, \frac{x}{\e_k}\right) = \bar{F}_k(M) = f_k(0)& \text{in }\Omega_{k,1},\\
v_k = - \e_k^2 w_{F_k} \left( M,\frac{x}{\e_k} \right) & \text{on }\partial\Omega_{k,1}. 
\end{dcases}
\end{equation}
Note that such a viscosity solution exists, since $\Omega_{k,1}$ satisfies a uniform exterior sphere condition, due to the assumption that $\Omega_k\in D(1,\bar\sigma)$. 

Since $\Omega_k\in D(1,\bar\sigma)$, we can assume, after extracting a subsequence if necessary, that $\Omega_k\to \Omega$ for some $\Omega\in D(1,\bar\sigma)$ in the sense of the Hausdorff distance. In particular, $\nu_k \to \nu$ for some $\nu\in \partial B_1$, and $\nu$ is the unit inward normal to $\partial\Omega$ at the origin. Also denoting by $\Phi_k$ and $\Phi$ the rotation mapping associated with $\Omega_k$ and respectively $\Omega$ (as in Definition \ref{definition:domain}) such that $\Phi_k(\nu_k) = e_n = \Phi(\nu)$, we have $\Phi_k \to \Phi$ in $\R^n$. 

Arguing similarly as in the proof of Lemma \ref{lemma:apprx-bdry-C1a}, one can also argue from the assumptions $G_k\in S_0(\lambda,\Lambda)$, $\Omega_k\in D(1,\bar\sigma)$, \eqref{eq:fk-Ca0-bdry}, \eqref{eq:gk-C2a0-bdry} and \eqref{eq:uk-pde-bdry} that there is $\bar{u} \in C_{loc}^\beta(\Omega_1\cup \Gamma_1)$ for which 
\begin{equation}\label{eq:uk-ub-bdry}
u_k \circ \Phi_k^{-1}\circ \Phi \to \bar{u}\quad\text{locally uniformly in }\Omega_1\cup \Gamma_1,\quad \frac{\partial u_k}{\partial \nu_k}(0) \to \frac{\partial \bar{u}}{\partial \nu} (0),
\end{equation}
possibly after extracting a subsequence. In particular, from \eqref{eq:gk-C2a0-bdry}, \eqref{eq:uk-pde-bdry} and the convergence above, we have 
\begin{equation}\label{eq:ub-bdry-1}
\bar{u} = 0 \quad\text{on }\Gamma_1,\quad |\bar{u}|\leq 1\quad\text{in }\Omega_1. 
\end{equation} 

In order to derive the interior equation for $\bar{u}$, we apply Lemma \ref{lemma:hom-stab2}. Due to \eqref{eq:fk-Ca0-bdry} and \eqref{eq:Gk-Fk-Linf-bdry}, we know that $f_k \circ \Phi_k^{-1} \circ \Phi \to 0$ uniformly in $\Omega_1$ and 
\begin{equation*}
\sup_{N\in\cS^n,y\in\R^n} \left(\frac{|G_k (N,\Phi_k^{-1}(\Phi(\cdot)),y) - F_k(N,y)|}{|N|}\right)  \to 0,\quad \text{locally uniformly in }\Omega_1.
\end{equation*} 
Thus, applying Lemma \ref{lemma:hom-stab2} to any fixed subdomain of $\Omega_1\cap\Omega_k$ and then letting $k\to\infty$, we deduce that $F_k \to F$ locally uniformly in $\cS^n\times\R^n$, for some periodic functional $F\in S_2(\lambda,\Lambda,\bar\kappa,\bar\gamma)\cap R_0(\kappa,\gamma)$ on $\cS^n\times\R^n$, 
\begin{equation}\label{eq:Fbk-Fb-bdry}
\bar{F}_k \to \bar{F} \quad\text{locally uniformly in }\cS^n,
\end{equation}
and that $\bar{u}$ is a viscosity solution of
\begin{equation}\label{eq:ub-pde-bdry11}
\bar{F}(D^2 \bar{u}) = 0  \quad\text{in }\Omega_1.
\end{equation}

In view of \eqref{eq:ub-bdry-1} and \eqref{eq:ub-pde-bdry11}, it follows from the standard boundary $C^{2,\bar\alpha}$ estimate \cite[Theorem 1.4]{SS} that there is $0<\bar\alpha<\bar\sigma$ and $\bar{c}>0$, depending only on $n$, $\lambda$, $\Lambda$ and $\bar\sigma$, such that $\bar{u} \in C^{2,\bar\alpha}(\overline{\Omega_{1/2}})$, and 
\begin{equation*}
\norm{\bar{u}}_{C^{2,\bar\alpha}(\overline{\Omega_{1/2}})} \leq \bar{c}.
\end{equation*} 
In particular, one can choose $0<\mu\leq\frac{1}{2}$, depending only on $n$, $\lambda$, $\Lambda$, $\bar{c}$, $\bar\alpha$ and $\alpha$, such that \begin{equation}\label{eq:ub-bdry-C2a-1}
\sup_{x\in \Omega_\mu} \left| \bar{u}(x) - \frac{\partial \bar{u}}{\partial \nu} (0) - \frac{1}{2} \inn{x}{\bar{M} x} \right| \leq \frac{1}{2} \mu^{2+\alpha},
\end{equation}
and moreover, (since $\nu_k \to \nu$)
\begin{equation}\label{eq:Mb-bdry-22}
|\bar{M}|\leq \bar{c},\quad (I - \nu\otimes\nu)\bar{M} = 0,\quad \bar{F}(\bar{M}) = 0.
\end{equation} 

Now let us construct a sequence $\{M_k\}_{k=1}^\infty\subset\cS^n$ such that each $M_k$ satisfies \eqref{eq:M-bdry-false} and $M_k\to \bar{M}$ as $k\to\infty$. As $\bar{F}_k$ being elliptic with ellipticity constants fixed by $\lambda$ and $\Lambda$ for all $k$'s, denoting by 
\begin{equation*}
\delta_k = |\bar{F}_k(\bar{M})|, \quad \rho_k =\lambda^{-1}( \delta_k + (\Lambda + 1)\e_k), 
\end{equation*} 
and noting that $\nu_k\otimes \nu_k \geq 0$ and $|\nu_k\otimes \nu_k| = 1$, one has 
\begin{equation*}
\bar{F}_k ( \bar{M} + D_T^2 g_k(0) + \rho_k (\nu_k\otimes \nu_k)) \geq \bar{F}_k (\bar{M}) - \Lambda |D_T^2 g_k(0)| + \delta_k + (\Lambda +1)\e_k \geq \e_k \geq f_k(0), 
\end{equation*}
due to \eqref{eq:fk-Ca0-bdry} and \eqref{eq:gk-C2a0-bdry}. Similarly, one obtains 
\begin{equation*}
\bar{F}_k (\bar{M} + D_T^2 g_k(0) - \rho_k (\nu_k\otimes \nu_k)) \leq - \e_k \leq f_k(0). 
\end{equation*}
Thus, the intermediate value theorem implies that for each $k=1,2,\cdots$, there exists some $0\leq t_k\leq 1$ such that  
\begin{equation}\label{eq:Mk-bdry11}
M_k = \bar{M} + D_T^2 g_k(0) + t_k \rho_k(\nu_k\otimes \nu_k),
\end{equation} 
satisfies
\begin{equation}\label{eq:Mk-bdry12}
\bar{F}_k (M_k) = f_k(0).
\end{equation}
Moreover, it follows from \eqref{eq:Mb-bdry-22} and \eqref{eq:Fbk-Fb-bdry} that 
\begin{equation}\label{eq:Mk-bdry13}
(I - \nu_k\otimes \nu_k) M_k = D_T^2 g_k(0),\quad M_k\to \bar{M}.
\end{equation} 
and that 
\begin{equation}\label{eq:Mk-bdry14}
|M_k|\leq \bar{C},
\end{equation}
for any large $k$, since $\bar{C}$ was given in the beginning of this proof by $\bar{C} = 2\bar{c} > \bar{c} \geq |M|$. 

Due to \eqref{eq:Mk-bdry14} and the fact that $F_k\in S_2(\lambda,\Lambda,\bar\kappa,\bar\gamma)\cap R_0(\kappa,\gamma)$, it follows from Lemma \ref{lemma:Fb} (iv) and \eqref{eq:Mk-bdry14} that 
\begin{equation}\label{eq:wFk-bdry-1}
\norm{w_{F_k}(M_k,\cdot)}_{C^{2,\gamma}(\R^n)} \leq C_0|M_k| \leq C_0\bar{C},
\end{equation} 
where $C_0$ depends at most on $n$, $\lambda$, $\Lambda$, $\bar\kappa$, $\bar\gamma$, $\kappa$, $\gamma$ and $\bar\sigma$. Now let $v_k$ be the viscosity solution of 
\begin{equation}\label{eq:vk-pde-bdry}
\begin{dcases}
G_k \left( M_k + D_y^2 w_{F_k}\left(M_k,\frac{x}{\e}\right) + D_x^2 v_k, x, \frac{x}{\e_k}\right) = \bar{F}_k(M_k) =f_k(0) & \text{in }\Omega_k,\\
v_k = - \e_k^2 w_{F_k} \left( M_k,\frac{x}{\e_k} \right) & \text{on }\partial\Omega_k.
\end{dcases}
\end{equation}
Note from \eqref{eq:wFk-bdry-1} that $v_k\in C^{0,1}(\partial\Omega_k)$ with
\begin{equation}\label{eq:vk-C01-bdry}
\norm{v_k}_{C^{0,1}(\partial\Omega_k)} \leq \e_k \norm{D_y w_{F_k}(M_k,\cdot)}_{L^\infty(\R^n)} \leq C_0\bar{C} \e_k. 
\end{equation}
Let us claim that 
\begin{equation}\label{eq:claim-vk-bdry}
v_k \circ \Phi_k^{-1} \circ \Phi \to 0 \quad\text{uniformly in }\Omega_1.
\end{equation}

If the claim \eqref{eq:claim-vk-bdry} is true, then from \eqref{eq:uk-ub-bdry}, \eqref{eq:Mk-bdry13}, \eqref{eq:wFk-bdry-1} and \eqref{eq:claim-vk-bdry}, we see that passing to the limit in \eqref{eq:apprx-bdry-false} with $M = M_k$ yields a contradiction to \eqref{eq:ub-bdry-C2a-1}. Thus, the proof will be finished, once we have the claim \eqref{eq:claim-vk-bdry}. 

To justify this claim, we shall construct a suitable barrier function. Utilizing the cell problem \eqref{eq:cell} associated with $F_k$ at $M_k$, we see that the interior equation in \eqref{eq:vk-pde-bdry} can be reformulated as 
\begin{equation}\label{eq:vk-pde-bdry-re}
H_k ( D^2 v_k,x) = h_k\quad\text{in }\Omega_{k,1},
\end{equation}
where $H_k$ and $h_k$ are defined by
\begin{equation*}
\begin{split}
H_k (N,x) &= G_k \left( M_k + D_y^2 w_{F_k} \left(M_k,\frac{x}{\e_k} \right) + N, x, \frac{x}{\e_k} \right) \\
&\quad - G_k \left(M_k + D_y^2 w_{F_k} \left(M_k,\frac{x}{\e_k}\right) ,x, \frac{x}{\e_k} \right)
\end{split}
\end{equation*}
and respectively
\begin{equation*} 
\begin{split}
h_k (x) &= G_k \left(M_k + D_y^2 w_{F_k} \left(M_k,\frac{x}{\e_k} \right),x, \frac{x}{\e_k} \right) \\
&\quad - F_k \left(M_k + D_y^2 w_{F_k} \left(M_k,\frac{x}{\e_k} \right), \frac{x}{\e_k} \right). 
\end{split}
\end{equation*}

Note that both $H_k$ and $h_k$ are well-defined, since $w_{F_k}(M_k,\cdot) \in C^{2,\gamma}(\R^n)$. Moreover, we have $H_k \in S_0(\lambda,\Lambda)$ on $\cS^n\times\Omega_{k,1}\times\R^n$, as $G_k \in S_0(\lambda,\Lambda)$. On the other hand, \eqref{eq:Gk-Fk-Lp-bdry} and \eqref{eq:wFk-bdry-1} together imply that
\begin{equation}\label{eq:hk-Lp-bdry}
\norm{h_k}_{L^p(\Omega_{k,1})} \leq \left( |M_k| + \norm{D_y^2 w_{F_k} (M_k,\cdot)}_{L^\infty(\R^n)} \right)L \leq (1+C_0)\bar{C}L,
\end{equation}
and similarly from \eqref{eq:Gk-Fk-Linf-bdry} with \eqref{eq:wFk-bdry-1} it follows that
\begin{equation}\label{eq:hk-Linf-bdry}
\norm{d(\cdot,\Gamma_{k,1})^2 h_k}_{L^\infty(\Omega_{k,1})} \leq (1 + C_0)\bar{C}\e_k.
\end{equation}

Now from \eqref{eq:vk-C01-bdry}, \eqref{eq:vk-pde-bdry-re} and \eqref{eq:hk-Lp-bdry}, we can invoke a global {\it a priori} estimate \cite[Proposition 4.14]{CC} such that $v_k\in C(\overline{\Omega_{k,1}})$ and there is a modulus of continuity $\rho$ such that 
\begin{equation}\label{eq:vk-cont-bdry}
|v_k(x) - v_k(x_0)| \leq \rho(|x-x_0|),\quad x,x_0\in\overline{\Omega_{k,1}}. 
\end{equation}
In particular, $\rho$ is determined only by $n$, the ellipticity constants of $H_k$, the diameter of $\Omega_{k,1}$, the radius with which $\Omega_{k,1}$ satisfies the uniform exterior sphere condition, the $L^n$ norm of $h_k$ on $\Omega_{k,1}$ and the modulus of continuity of $v_k$ on $\partial\Omega_{k,1}$. Thus, the dependence of $\rho$ reduces to the parameters $n$, $\lambda$, $\Lambda$, $\bar\sigma$, $C_0$, $\bar{C}$ and $L$; especially, it is independent of $k$, and it also has nothing do with either $\nu_k$ and $\nu$ being rational or irrational direction.

Let $0<\delta<1$ be arbitrary. Then from \eqref{eq:vk-cont-bdry} and \eqref{eq:vk-C01-bdry} we know that 
\begin{equation}\label{eq:claim-vk-bdry-1}
|v_k(x)| \leq \rho(\delta) + C_0 \bar{C} \e_k,\quad \text{if } d(x,\Gamma_{k,1}) \leq \delta. 
\end{equation}
On the set $\Omega_{k,1} \cap \{d(\cdot,\Gamma_{k,1})\geq \delta\}$, consider an auxiliary function
\begin{equation*}
\psi_k (x) = \frac{(1 + C_0)\bar{C}\e_k}{2\lambda\delta^2} (1 - |x|^2) + \rho(\delta) + C_0 \bar{C} \e_k. 
\end{equation*}
Clearly, on the boundary, we deduce from  \eqref{eq:claim-vk-bdry-1} that
\begin{equation}\label{eq:claim-vk-bdry-2}
\psi_k (x) \geq v_k (x),\quad \text{if } d(x,\Gamma_{k,1}) = \delta.
\end{equation} 
On the other hand, in the interior, from the fact that $H_k\in S_0(\lambda,\Lambda)$ as well as the estimate \eqref{eq:hk-Linf-bdry} for $h_k$, it follows that 
\begin{equation}\label{eq:claim-vk-bdry-3}
H_k (D^2 \psi_k(x),x) \leq - \frac{(1+C_0)\bar{C} \e_k}{\delta^2} \leq h_k(x),\quad \text{if }d(x,\Gamma_{k,1}) >\delta. 
\end{equation}
Hence, $\psi_k$ is a supersolution to the boundary value problem \eqref{eq:vk-pde-bdry-re} that $v_k$ solves in the viscosity sense, so it follows from the comparison principle that
\begin{equation}\label{eq:claim-vk-bdry-4}
v_k (x) \leq \psi_k (x) \leq \left( C_0 + \frac{1+C_0}{2\lambda\delta^2} \right) \bar{C} \e_k + \rho(\delta),\quad \text{if } d(x,\Gamma_{k,1}) \geq \delta. 
\end{equation} 
Letting $k\to\infty$ in both \eqref{eq:claim-vk-bdry-1} and \eqref{eq:claim-vk-bdry-4}, with $\delta$ fixed, we obtain 
\begin{equation*}
\limsup_{k\to\infty} v_k\circ\Phi_k^{-1}\circ\Phi \leq \rho(\delta)\quad\text{in }\Omega_1. 
\end{equation*}
Since $\delta$ can be arbitrarily small and $\rho$ is a modulus of continuity, the last inequality implies that $\limsup_{k\to\infty} v_k\circ\Phi_k^{-1}\circ\Phi\leq 0$ in $\Omega_1$. By a similar argument, one can also show that $\liminf_{k\to\infty} v_k\circ\Phi_k^{-1}\circ\Phi \geq 0$ in $\Omega_1$, verifying the claim \eqref{eq:claim-vk-bdry}. 
\end{proof}

\begin{remark}\label{remark:apprx-bdry}
If $L = 0$, hence $p = \infty$ in Lemma \ref{lemma:apprx-bdry}, $G$ coincides with $F$ everywhere. Then we can replace Lemma \ref{lemma:hom-stab2} with Lemma \ref{lemma:hom-stab}. Since the latter works with the class $S_0(\lambda,\Lambda)$, the dependence of $\bar\e$ can be restricted to the parameters $n$, $\lambda$, $\Lambda$, $\bar\sigma$ and $\alpha$. This will be used in the initial step ($k = 0$) in the iteration lemma below. From the second step ($k\geq 1$), one has to work with $L>0$ and some finite $p>n$. 
\end{remark}

Next follows an iteration lemma. 

\begin{lemma}\label{lemma:iter-bdry}
Let $\bar\sigma$, $\bar\alpha$, $\bar{C}$, $\alpha$ and $\mu$ be as in Lemma \ref{lemma:apprx-bdry}. Let $F \in S_2(\lambda,\Lambda,\bar\kappa,\bar\gamma)\cap R_1(\kappa,\gamma)$ be a periodic functional on $\cS^n\times\R^n$, and suppose that $\Omega\in D(\eta,\bar\sigma)$, $f\in C^\alpha(\overline{\Omega_1})$, $g\in C^{2,\alpha}(\Gamma_1)$, and $u^\e\in C(\overline{\Omega_1})$ satisfy
\begin{equation}\label{eq:g-C11-bdry}
g(0) = |D_T g(0)| = 0,
\end{equation} 
\begin{equation}\label{eq:Dnue-0}
\frac{\partial u^\e}{\partial\nu} (0) = 0,\quad\text{and}
\end{equation}
\begin{equation}\label{eq:ue-pde-bdry-gen}
\begin{dcases}
F \left( D^2 u^\e, \frac{x}{\e}\right) = f & \text{in }\Omega_1,\\
u^\e = g & \text{on }\Gamma_1.
\end{dcases}
\end{equation}
Also let $K>0$ and $p>n$ be given. Then there are $0<\eta,\bar\e\leq\frac{1}{2}$, depending only on $n$, $\lambda$, $\Lambda$, $\bar\sigma$ and $\alpha$, and $0<\hat\e\leq \bar\e$, depending only on $n$, $\lambda$, $\Lambda$, $\bar\kappa$, $\bar\gamma$, $\kappa$, $\gamma$, $\bar\sigma$, $\alpha$, $K$ and $p$, such that if $\e\leq \hat\e \mu^{k-1}$ for some integer $k\geq 1$, $\Omega\in D(\eta,\bar\sigma)$ and 
\begin{equation}\label{eq:Ke-bdry}
J^\e = \norm{u^\e}_{L^\infty(\Omega_1)} + \frac{1}{\bar\e } \norm{ f }_{C^\alpha(\overline{\Omega_1})} + \frac{4}{\bar\e } \norm{g}_{C^{2,\alpha}(\Gamma_1)} \leq K,
\end{equation}
then there are $a_k^\e \in \R$ and $M_k^\e\in \cS^n$, satisfying
\begin{equation}\label{eq:ake-Mke-bdry}
|a_k^\e| \leq \frac{C_1\bar{C}\e^{1-\beta}}{1-\mu^\alpha} J^\e,\quad |M_k^\e| \leq \frac{\bar{C}}{1-\mu^\alpha} J^\e ,\quad \bar{F}(M_k^\e) = f(0), 
\end{equation}
as well as a viscosity solution $\zeta_k^\e$ to
\begin{equation}\label{eq:zke-pde-bdry}
\begin{dcases}
F \left( M_k^\e + D_y^2 w_F \left(M_k^\e,\frac{x}{\e}\right) + D_x^2 \zeta_k^\e,\frac{x}{\e}\right) = \bar{F}(M_k^\e) = f(0) & \text{in }\Omega_{\mu^{k-1}}, \\
\zeta_k^\e = -\e^2 w_F \left( M_k^\e,\frac{x}{\e}\right) & \text{on }\Gamma_{\mu^{k-1}},\\
|\zeta_k^\e| \leq \frac{C_2\bar{C}\e^2}{1-\mu^\alpha}J^\e & \text{on }\partial\Omega_{\mu^k}\setminus\Gamma_{\mu^{k-1}},
\end{dcases}
\end{equation}
such that 
\begin{equation}\label{eq:iter-bdry}
\sup_{x\in \Omega_{\mu^k}} \left| u^\e(x) - a_k^\e \inn{x}{\nu} - \frac{1}{2} \inn{x}{M_k^\e x} - \e^2 w_F \left( M_k^\e,\frac{x}{\e}\right) - \zeta_k^\e(x) \right| \leq J^\e \mu^{k(2+\alpha)},
\end{equation}
where $\beta<1$, $C_1$ and $C_2$ depend only on $n$, $\lambda$ and $\Lambda$.
\end{lemma}

\begin{proof}
Let us choose $\bar\e$ such that Lemma \ref{lemma:apprx-bdry} holds with $L = 0$. As noted in Remark \ref{remark:apprx-bdry}, such $\bar\e$ depends only on $n$, $\lambda$, $\Lambda$, $\bar\sigma$ and $\alpha$. Let $\hat\e$, $\eta$, $\beta$ and $C_1$ be determined later. 

Henceforth, suppose that $\e\leq \hat\e \mu^k$ for some integer $k\geq 0$, and we have already found some $a_k^\e\in\R$, $M_k^\e\in\cS^n$ and a harmonic function $\zeta_k^\e$ in $\Omega_{\mu^k}$ satisfying \eqref{eq:iter-bdry}. If $k=0$, we choose $a_0^\e=0$, $M_0^\e = 0$ and $\zeta_0^\e =0$, so that \eqref{eq:iter-bdry} clearly holds. If $k\geq 1$, let us assume that $a_k^\e\in\R$ satisfies the first inequality in \eqref{eq:ake-Mke-bdry}, and that $M_k^\e\in\cS^n$ satisfies 
\begin{equation}\label{eq:ake-bdry-re}
|(I - \nu\otimes \nu) M_k^\e  - D_T^2 g(0)| \leq \frac{C_1\bar{C}\eta\e^{1-\beta}}{1-\mu^\alpha}J^\e,\quad |M_k^\e |\leq \bar{C} J^\e \sum_{l=1}^k \mu^{(l-1)\alpha},
\end{equation}
as well as the last equality in \eqref{eq:ake-Mke-bdry}. Moreover, suppose that $\zeta_k^\e$ the viscosity solution to \eqref{eq:zke-pde-bdry}. 

Define, for $x\in \mu^{-k} \Omega_{\mu^k} = \{ x\in B_1: \mu^k x\in \Omega\}$, 
\begin{equation}\label{eq:uke-bdry}
u_k^\e (x) = \frac{u^\e(\mu^k x) - \mu^k a_k^\e \inn{x}{\nu} - \frac{\mu^{2k}}{2} \inn{x}{M_k^\e x} - \e^2 w_F (M_k^\e,\e^{-1}\mu^k x) - \zeta_k^\e(\mu^k x)}{J^\e \mu^{k(2+\alpha)}}
\end{equation}
and 
\begin{equation}\label{eq:fke-Vke}
f_k^\e (x) = \frac{f(\mu^k x ) - f(0)}{J^\e \mu^{k\alpha}}.
\end{equation}
Also set, for $x\in \mu^{-k}\Gamma_{\mu^k} = \{ x\in B_1: \mu^k x\in\Gamma\}$,  
\begin{equation}\label{eq:gke-bdry}
g_k^\e (x)  = \frac{g(\mu^k x) - \mu^k a_k^\e \inn{x}{\nu} - \frac{\mu^{2k}}{2} \inn{x}{M_k^\e x}}{J^\e \mu^{k(2+\alpha)}}.
\end{equation}
As $u^\e$ being a viscosity solution of \eqref{eq:ue-pde-bdry-gen}, and since $u^\e$ satisfies \eqref{eq:iter-bdry} as an induction hypothesis, we have 
\begin{equation}\label{eq:uke-pde-bdry}
\begin{dcases}
G_k^\e \left( D^2 u_k^\e,x,\frac{\mu^k x}{\e}\right) = f_k^\e,\quad |u_k^\e| \leq 1 & \text{in }\mu^{-k}\Omega_{\mu^k},\\
u_k^\e = g_k^\e & \text{on }\mu^{-k} \Gamma_{\mu^k},
\end{dcases}
\end{equation}
in the viscosity sense, where $G_k^\e$ is defined by
\begin{equation}\label{eq:Gke-bdry}
\begin{split}
G_k^\e (N,x,y) & = \frac{1}{J^\e \mu^{k\alpha}} F (J^\e \mu^{k\alpha} N + M_k^\e + D_y^2 w_F (M_k^\e,y) + D_x^2 \zeta_k^\e (\mu^k x), y) \\
&\quad - \frac{1}{J^\e \mu^{k\alpha}} F (M_k^\e + D_y^2 w_F (M_k^\e,y) + D_x^2 \zeta_k^\e (\mu^k x),y),
\end{split}
\end{equation}
so that $G_k^\e$ is a periodic functional on $\cS^n\times(\mu^{-k}\Omega_{\mu^k})\times\R^n$ belonging to $S_0(\lambda,\Lambda)$.

From \eqref{eq:fke-Vke} and \eqref{eq:Ke-bdry}, it is clear that
\begin{equation}\label{eq:fke-Ca0-bdry}
\norm{f_k^\e}_{C^\alpha(\mu^{-k}\overline{\Omega_{\mu^k}})} \leq \frac{1}{J^\e} [f]_{C^\alpha(\overline{\Omega_{\mu^k}})} \leq \bar\e .
\end{equation}
On the other hand, we also have 
\begin{equation}\label{eq:gke-C2a0-bdry}
\norm{g_k^\e}_{C^{2,\alpha}(\mu^{-k} \Gamma_{\mu^k})} \leq \bar\e .
\end{equation}
Let us stress that this estimate is irrelevant to the nonlinear structure in the interior homogenization, and moreover it has nothing to do with the PDE that $\zeta_k^\e$ satisfies in the interior. Hence, one can follow exactly the same argument in \cite[Lemma 11]{AL2}. Since the argument is long and technical, we shall not repeat it here. Still let us remark that a direct computation yields 
\begin{equation}\label{eq:DTgke-bdry}
g_k^\e (0) = | D_T g_k^\e (0) | = 0,\quad\text{and}
\end{equation}
\begin{equation}\label{eq:DT2gke-bdry}
D_T^2 g_k^\e (0) = -\frac{1}{J^\e \mu^{k\alpha}} (a_k^\e D_T^2 \phi(0) + (I - \nu\otimes\nu)M_k^\e - D_T^2 g(0)), 
\end{equation}
where $\phi:\Pi \to \R$ is the parameterization of $\Gamma_1$ with respect to the hyperplane $\Pi_1 = \{x\in B_1:\inn{x}{\nu}=0\}$ such that 
\begin{equation*}
\Gamma_1 = \{z+ \phi(z)\nu: z\in \Pi_1\},\quad\text{and}
\end{equation*}
\begin{equation*}
\phi(0) = |D_T \phi(0)| = 0,\quad \norm{\phi}_{C^{2,\bar\sigma}(\Pi\cap B_1)} \leq \eta. 
\end{equation*}
Such a characterization of $\Gamma_1$ exists uniquely up to a rotation which fixes the direction $\nu$, due to the assumption $\Omega\in D(\eta,\bar\sigma)$. Following the argument in \cite[Lemma 11]{AL2} carefully, one can also observe that the smallness condition of $\eta\leq 1$ is determined only by $C_1$, $\bar{C}$, $\mu$, $\alpha$ and $\bar\e$, hence on $n$, $\lambda$, $\Lambda$, $\bar\sigma$ and $\alpha$, provided that $C_1$ depends at most on $n$, $\lambda$ and $\Lambda$; we shall choose $C_1$ at the end of this proof. 

In order to verify that one can apply the approximation lemma, Lemma \ref{lemma:apprx-bdry}, let us define $F_k^\e$ on $\cS^n\times\R^n$ by 
\begin{equation}\label{eq:Fke-bdry}
F_k^\e (N,y) = \frac{F(J^\e \mu^{k\alpha} N + M_k^\e + D_y^2 w_F (M_k^\e,y),y) - f(0)}{J^\e \mu^{k\alpha}}.
\end{equation}
From the induction hypothesis $\bar{F}(M_k^\e) = f(0)$, we also have the alternative definition \eqref{eq:Fke-re}. Hence, $F_k^\e$ is a periodic functional on $\cS^n\times\R^n$ belonging to $S_0(\lambda,\Lambda)$. Moreover, since $F\in S_2(\lambda,\Lambda,\bar\kappa,\bar\gamma)\cap R_1(\kappa,\gamma)$, it follows from the second assertion in Lemma \ref{lemma:GF} (or as in the proof of Theorem \ref{theorem:int-C1a} (iii)) that
\begin{equation}\label{eq:Fke-class-bdry}
F_k^\e \in R_0 (\kappa(C_0|M_k^\e| + 1),\gamma) \subset R_0 \left( \kappa\left(\frac{C_0\bar{C}K}{1-\mu^\alpha}  +1 \right),\gamma \right),
\end{equation}
with $C_0$ depending only on $n$, $\lambda$, $\Lambda$, $\bar\kappa$, $\bar\gamma$, $\kappa$ and $\gamma$, where the second inclusion is due to the induction hypothesis \eqref{eq:ake-bdry-re} on $M_k^\e$ and the assumption \eqref{eq:Ke-bdry} on $J^\e$.

On the other hand, in view of \eqref{eq:Fke-re}, the interior equation in \eqref{eq:zke-pde-bdry} of $\zeta_k^\e$ can be reformulated as 
\begin{equation}\label{eq:zke-pde-bdry-re}
F_k^\e \left( D^2 \zeta_k^\e,\frac{x}{\e}\right) = 0 \quad\text{in }\Omega_{\mu^{k-1}}.
\end{equation}
This together with the boundary condition in \eqref{eq:zke-pde-bdry} yields from the Alexandroff-Bakelman-Pucci estimate \cite[Theorem 3.6]{CC} and the $L^\infty$ estimate \eqref{eq:wF-Ca} of $w_F(M_k^\e,\cdot)$ that 
\begin{equation}\label{eq:zke-Linf-bdry}
\norm{\zeta_k^\e}_{L^\infty(\Omega_{\mu^{k-1}})} \leq \frac{C_2\bar{C}\e^2}{1-\mu^\alpha}J^\e, 
\end{equation}
provided that $C_2$ is chosen by the constant in \eqref{eq:wF-Ca}, which depends only on $n$, $\lambda$ and $\Lambda$. Moreover, it follows from the {\it a priori} gradient estimate, such as \cite[Theorem 1.1]{SS}, and the $C^{1,\beta}$ estimate \eqref{eq:wF-C1a} of $w_F(M_k^\e,\cdot)$ that 
\begin{equation}\label{eq:zke-Lip-bdry}
\begin{split}
\left| \frac{\partial \zeta_k^\e}{\partial \nu}(0) \right| &\leq C_3\left\{ \norm{\zeta_k^\e}_{L^\infty(\Omega_{\mu^{k-1}})} + \e^2 \norm{w_F \left(M_k^\e,\frac{\cdot}{\e}\right)}_{C^{1,\beta}(\Gamma_{\mu^{k-1}})} \right\} \leq \frac{C_4\bar{C}\e^{1-\beta}}{1-\mu^\alpha} J^\e,
\end{split}
\end{equation} 
where $\beta<1$, $C_3$ and $C_4$, depending only on $n$, $\lambda$ and $\Lambda$; here we used the assumption that $\Omega \in D(\eta,\bar\sigma)\subset D(1,\bar\sigma)$, implying that the maximal curvature of $\Gamma$ is bounded by $1$. 

In addition, in view of \eqref{eq:zke-pde-bdry-re} and the boundary condition in \eqref{eq:zke-pde-bdry} for $\zeta_k^\e$, one can apply the boundary $W^{2,p}$ estimate \eqref{eq:bdry-W2p} with $p>n$, after a scaling argument, and deduce from \eqref{eq:zke-Linf-bdry} and the $C^{2,\gamma}$ estimate \eqref{eq:wF-C2a} of $w_F (M_k^\e,\cdot)$ that
\begin{equation*}
\begin{split}
\norm{D^2 \zeta_k^\e}_{L^p(\Omega_{\mu^{k-1}/2})}& \leq C_5 \mu^{(k-1)n/p} \left\{ \frac{\norm{\zeta_k^\e}_{L^\infty(\Omega_{\mu^{k-1}})}}{\mu^{2(k-1)}} + \e^2 \norm{w_F\left(M_k^\e,\frac{\cdot}{\e}\right)}_{C^{1,1}(\Gamma_{\mu^{k-1}})} \right\} \\
&\leq C_6 \mu^{(k-1)n/p}\left( \frac{\bar{C}\e^2 J^\e}{(1-\mu^\alpha)\mu^{2(k-1)}} + |M_k^\e|\right) \\
&\leq \frac{C_7\bar{C} \mu^{(k-1)n/p+2}}{1-\mu^\alpha} J^\e, 
\end{split}
\end{equation*} 
with $C_5$, $C_6$ and $C_7$ depending at most on $n$, $\lambda$, $\Lambda$, $\bar\kappa$, $\bar\gamma$, $\kappa$, $\gamma$, $K$, $C_0$ and $C_2$, hence on the first eight parameters only, where in the third inequality we used the assumption $\e \mu^k \leq \hat\e\leq\frac{1}{2}$. In particular, since $\mu\leq\frac{1}{2}$, we have $\Omega_{\mu^k} \subset \Omega_{\mu^{k-1}/2}$ and 
\begin{equation}\label{eq:zke-W2p-bdry}
\norm{D^2 \zeta_k^\e}_{L^p(\Omega_{\mu^k})} \leq \frac{C_7 \bar{C} \mu^{kn/p}}{1-\mu^\alpha}J^\e.
\end{equation} 

On the other hand, we also have the uniform interior $C^{1,1}$ estimate \eqref{eq:int-C11} for $\zeta_k^\e$. This along with \eqref{eq:zke-Linf-bdry} one can deduce that
\begin{equation*}
\norm{d(\cdot,\partial\Omega_{\mu^{k-1}})^2 D^2 \zeta_k^\e}_{L^\infty(\Omega_{\mu^{k-1}})} \leq \frac{C_8\bar{C}\e^2}{1-\mu^\alpha} J^\e,
\end{equation*} 
where $C_8$ depends at most on $n$, $\lambda$, $\Lambda$, $\bar\kappa$, $\bar\gamma$, $\kappa$ and $\gamma$. Again since $\Omega_{\mu^k}\subset \Omega_{\mu^{k-1}/2}$, we have $d(\cdot,\partial\Omega_{\mu^{k-1}}) = d(\cdot,\Gamma_{\mu^k})$ in $\Omega_{\mu^k}$, which implies that 
\begin{equation}\label{eq:zke-C11-bdry}
\norm{d(\cdot,\Gamma_{\mu^k})^2 D^2 \zeta_k^\e}_{L^\infty(\Omega_{\mu^k})} \leq \frac{C_8\bar{C}\e^2}{1-\mu^\alpha} J^\e. 
\end{equation} 

Now we are ready to verify that $F_k^\e$ and $G_k^\e$ also satisfy \eqref{eq:G-F-Lp-bdry} and \eqref{eq:G-F-Linf-bdry} with some $\hat{L}$ possibly larger than $L$, such that we can apply the approximation lemma, Lemma \ref{lemma:apprx-bdry}. Note that one can write $F_k^\e$ and $G_k^\e$ by
\begin{equation*}
F_k^\e (N,y) = \tr(A_k^\e (N,y) N),\quad G_k^\e (N,x,y) = \tr(B_k^\e (N,x,y) N),
\end{equation*} 
where $A_k^\e$ and $B_k^\e$ are defined by
\begin{equation*}
A_k^\e (N,y) = \int_0^1 D_M F (t J^\e \mu^{k\alpha} N + M_k^\e + D_y^2 w_F (M_k^\e,y),y) \,dt,
\end{equation*} 
and respectively by
\begin{equation*}
B_k^\e (N,x,y) = \int_0^1 D_M F (t J^\e \mu^{k\alpha} N + M_k^\e + D_y^2 w_F (M_k^\e,y) + D_x^2 \zeta_k^\e(\mu^k x),y) \,dt.
\end{equation*} 
Thus, the structure condition \eqref{eq:F-C11} on $F$ implies that 
\begin{equation*}
| G_k^\e (N,x,y) - F_k^\e (N,y) | \leq | B_k^\e (N,x,y) - A_k^\e (N,y) | |N| \leq \kappa|D^2 \zeta_k^\e (\mu^k x)||N|,
\end{equation*}
for any $N\in\cS^n$, $x\in\mu^{-k}\Omega_{\mu^k}$ and $y\in\R^n$. This estimate combined with \eqref{eq:zke-W2p-bdry} and \eqref{eq:zke-C11-bdry} yields that
\begin{equation}\label{eq:Gke-Fke-Lp-bdry}
\norm{G_k^\e (N,\cdot,y) - F_k^\e (N,y)}_{L^p(\mu^{-k}\Omega_{\mu^k})} \leq \left(\frac{\kappa C_7\bar{C}\mu^{2-n/p}J^\e}{1-\mu^\alpha} \right)|N|,
\end{equation}
and respectively 
\begin{equation}\label{eq:Gke-Fke-Linf-bdry}
\norm{d(\cdot,\partial (\mu^{-k}\Gamma_{\mu^k}))^2 (G_k^\e (N,\cdot,y) - F_k^\e (N,y))}_{L^\infty(\mu^{-k}\Omega_{\mu^k})} \leq \left( \frac{\kappa C_8\bar{C}J^\e}{1-\mu^\alpha}\right)\left(\frac{\e}{\mu^k}\right)^2 |N| .
\end{equation}

We are finally in a position to determine the smallness condition on $\hat\e$. Let us first take $\hat\e\leq \bar\e$ in such a way that Lemma \ref{lemma:apprx-bdry} holds with  
\begin{equation}\label{eq:Lh-bdry}
L = \frac{\kappa C_7\bar{C}\mu^{2-n/p}K}{1-\mu^\alpha}
\end{equation} 
and the parameters involved in the class $R_0$ for $F_k^\e$, given by \eqref{eq:Fke-class-bdry}, as well as $\bar\sigma$, $\alpha$ and $p$. Then $\hat\e$ depends only on $n$, $\lambda$, $\Lambda$, $\bar\kappa$, $\bar\gamma$, $\kappa$, $\gamma$, $\bar\sigma$, $\alpha$, $K$ and $p$. Then we let $\hat\e$ even smaller, if necessary,  so as to satisfy
\begin{equation}\label{eq:eh-bdry}
\left( \frac{\kappa C_8\bar{C}J^\e}{1-\mu^\alpha}\right)\left(\frac{\e}{\mu^k}\right)^2 \leq \left( \frac{\kappa C_8\bar{C}K}{1-\mu^\alpha}\right)\hat\e^2 \leq \bar\e,
\end{equation} 
where the first inequality follows from the assumption that $\e\leq \hat\e \mu^k$. This will not change the dependence of $\hat\e$ specified above. 

From \eqref{eq:fke-Ca0-bdry}, \eqref{eq:gke-C2a0-bdry}, \eqref{eq:Gke-Fke-Lp-bdry}, \eqref{eq:Gke-Fke-Linf-bdry}, \eqref{eq:Lh-bdry} and \eqref{eq:eh-bdry}, as well as the fact that $\mu^{-k}\Omega_{\mu^k} \in D(\eta,\bar\sigma)\subset D(1,\bar\sigma)$, we can apply Lemma \ref{lemma:apprx-bdry} to the problem \eqref{eq:uke-pde-bdry}, with $\e$ replaced by $\e\mu^{-k}$. This yields $N_k^\e\in\cS^n$, which satisfies
\begin{equation}\label{eq:Nke-bdry}
|N_k^\e| \leq \bar{C}, \quad (I - \nu\otimes \nu) N_k^\e = D_T^2 g_k^\e(0), \quad \bar{F}_k^\e (N_k^\e) = f_k^\e(0) = 0,
\end{equation}
and the viscosity solution $v_k^\e$ to 
\begin{equation}\label{eq:vke-pde-bdry}
\begin{dcases}
G_k^\e \left( N_k^\e + D_y^2 w_{F_k^\e} \left(N_k^\e,\frac{\mu^k x}{\e}\right) + D_x^2 v_k^\e, x, \frac{\mu^k x}{\e}\right) = \bar{F}_k^\e (N_k^\e ) = 0 & \text{in }\mu^{-k}\Omega_{\mu^k},  \\
v_k^\e = -\frac{\e^2}{\mu^{2k}} w_{F_k^\e} \left(N_k^\e,\frac{\mu^k x}{\e}\right) & \text{on }\partial(\mu^{-k}\Omega_{\mu^k}),
\end{dcases}
\end{equation}
such that for any $x\in\mu^{-k}\Omega_{\mu^k}\cap B_\mu = \mu^{-k} \Omega_{\mu^{k+1}}$, 
\begin{equation}\label{eq:uke-iter-1}
\left| u_k^\e (x) - \frac{\partial u_k^\e}{\partial\nu}(0) \inn{x}{\nu} - \frac{1}{2} \inn{x}{N_k^\e x} - \frac{\e^2}{\mu^{2k}} w_{F_k^\e}\left( N_k^\e,\frac{\mu^k x}{\e}\right) - v_k^\e (x)\right| \leq \mu^{2+\alpha}.
\end{equation}
Note from the Alexandroff-Bakelman-Pucci estimate \cite[Theorem 3.6]{CC} that 
\begin{equation}\label{eq:vke-Linf-bdry}
\norm{v_k^\e}_{L^\infty(\mu^{-k}\Omega_{\mu^k})} \leq \frac{\e^2}{\mu^{2k}} \norm{w_{F_k^\e}\left(N_k^\e,\cdot\right)}_{L^\infty(\partial \Omega_{\mu^k})} \leq \frac{C_2\bar{C}\e^2}{\mu^{2k}},
\end{equation}
where we used the fact that $G_k^\e \in S_0(\lambda,\Lambda)$, the first inequality in \eqref{eq:Nke-bdry} and the $L^\infty$ estimate \eqref{eq:wF-Ca} for $w_{F_k^\e}(N_k^\e,\cdot)$; in particular, we chose $C_2$ by the constant appearing in \eqref{eq:wF-Ca}, which depends only on $n$, $\lambda$ and $\Lambda$. 

To this end, we define 
\begin{equation}\label{eq:ake-bdry-re2}
a_{k+1}^\e = a_k^\e + J^\e \mu^{k(1+\alpha)} \frac{\partial u_k^\e}{\partial \nu} (0) =  - \e \nu\cdot D_y w_F(M_k^\e,0) - \frac{\partial  \zeta_k^\e}{\partial \nu} (0),
\end{equation}
and
\begin{equation}\label{eq:Mke-bdry-re3}
M_{k+1}^\e = M_k^\e + J^\e \mu^{k\alpha} N_k^\e.
\end{equation}
Let us also define $\zeta_{k+1}^\e$ on $\overline{\Omega_{\mu^k}}$ by 
\begin{equation}\label{eq:zke-bdry-re3}
\zeta_{k+1}^\e (x) = \zeta_k^\e (x) + J^\e \mu^{k(2+\alpha)} v_k^\e \left( \frac{x}{\mu^k} \right).
\end{equation}
Due to \eqref{eq:Mke-bdry-re3}, the definition \eqref{eq:Fke-bdry} of $F_k^\e$ and Lemma \ref{lemma:eff-cor}, we have the additive structure \eqref{eq:Gbke} of $\bar{F}_k^\e$ and $w_{F_k^\e}$, as in the proof of Lemma \ref{lemma:iter}. Especially, we have 
\begin{equation}\label{eq:Fbke-bdry}
\bar{F}(M_{k+1}^\e) = \bar{F}(M_k^\e) = f(0), \quad w_F (M_{k+1}^\e,y) = w_F (M_k^\e ,y) + J^\e \mu^{k\alpha} w_{F_k^\e}(N_k^\e,y),
\end{equation}
where in the first equality we used \eqref{eq:Nke-bdry} and \eqref{eq:ake-Mke-bdry}. Hence, one can rephrase \eqref{eq:uke-iter-1} in terms of $u^\e$, as we have, for all $x\in \Omega_{\mu^{k+1}}$, 
\begin{equation*}
\left| u^\e (x) - a_{k+1}^\e \inn{x}{\nu} - \frac{1}{2} \inn{x}{M_{k+1}^\e x} - \e^2 w_F \left(M_{k+1}^\e,\frac{x}{\e}\right) - \zeta_{k+1}^\e (x) \right| \leq J^\e \mu^{(k+1)(2+\alpha)}.
\end{equation*}
This estimate verifies the induction hypothesis \eqref{eq:iter-bdry} with newly obtained $a_{k+1}^\e$, $M_{k+1}^\e$ and $\zeta_{k+1}^\e$. Thus, the proof is finished if one verifies \eqref{eq:ake-Mke-bdry}, \eqref{eq:ake-bdry-re} and \eqref{eq:zke-pde-bdry} for $k+1$. 

The first inequality in \eqref{eq:ake-Mke-bdry} for $k+1$ follows immediately from \eqref{eq:zke-Lip-bdry} and \eqref{eq:wF-C1a}, if we select $C_1$ by $2C_4$, with $C_4$ as in \eqref{eq:zke-Lip-bdry}; hence, $C_1$ depends only on $n$, $\lambda$ and $\Lambda$. This shows that $a_{k+1}^\e$ verifies its induction hypothesis. 

Regarding $M_{k+1}^\e$, it is clear that the second inequality in \eqref{eq:ake-bdry-re} for $k+1$ holds, owing to the induction hypothesis \eqref{eq:ake-bdry-re} for $M_k^\e$, and the first inequality in \eqref{eq:Nke-bdry} for $N_k^\e$. On the other hand, the first inequality in \eqref{eq:ake-bdry-re} for $k+1$ can be deduced from the second identity in \eqref{eq:Nke-bdry}, the first inequality in \eqref{eq:ake-Mke-bdry} and the observation \eqref{eq:DT2gke-bdry}. The last identity in \eqref{eq:ake-Mke-bdry} for $k+1$ is already verified by \eqref{eq:Fbke-bdry}. Therefore, $M_{k+1}^\e$ also satisfies its induction hypotheses.

Finally, the proof will be finished if we verify that $\zeta_{k+1}^\e$ solves to the boundary value problem \eqref{eq:zke-pde-bdry} for $k+1$. Utilizing \eqref{eq:Mke-bdry-re3}, \eqref{eq:Fbke-bdry}, the interior equations in \eqref{eq:vke-pde-bdry} and \eqref{eq:zke-pde-bdry} that $v_k^\e$ and $\zeta_k^\e$ satisfy respectively, one can proceed as
\begin{equation*}
\begin{split}
& F \left( M_{k+1}^\e + D_y^2 w_F \left( M_{k+1}^\e,\frac{x}{\e}\right) + D_x^2 \zeta_{k+1}^\e (x),\frac{x}{\e}\right) \\
& = F \left( M_k^\e + D_y^2 w_F \left( M_k^\e,\frac{x}{\e}\right) + D_x^2 \zeta_k^\e (x)  \right. \\
&\quad \quad \quad \left. + J^\e \mu^{k\alpha} \left( N_k^\e + D_y^2 w_{F_k^\e} \left( N_k^\e,\frac{x}{\e}\right) + D_x^2 v_k^\e \left(\frac{x}{\mu^k}\right)\right),\frac{x}{\e}\right) \\
& = F \left( M_k^\e + D_y^2 w_F \left(M_k^\e,\frac{x}{\e}\right) + D_x^2 \zeta_k^\e (x),\frac{x}{\e}\right) \\
& = f(0),
\end{split}
\end{equation*}
for $x\in \Omega_{\mu^k}$, which verifies that $\zeta_{k+1}^\e$ satisfies the interior equation in \eqref{eq:zke-pde-bdry} for $k+1$. The boundary condition for $\zeta_{k+1}^\e$ can be verified in a similar way. Thus, the proof is finished. 
\end{proof}

We are ready to prove the uniform boundary $C^{1,1}$ estimates. 

\begin{proof}[Proof of Theorem \ref{theorem:bdry-C11}]
With the iteration lemma (Lemma \ref{lemma:iter-bdry}) at hand, the proof is similar to \cite[Theorem 1]{AL2}, whose argument can be easily extended to fully nonlinear equations. Some necessary detail adopting the nonlinear structure can also be found in the proof of Theorem \ref{theorem:int-C11}. For this reason, we shall omit the detail and finish the proof here.
\end{proof}


\section{Examples}\label{section:ex}

In this section, we shall present some classes of periodically oscillating fully nonlinear functionals $F$ that verify the assumptions for the uniform $C^{1,1}$ estimates, namely Theorem \ref{theorem:int-C11} and Theorem \ref{theorem:bdry-C11}. The key assumption for these theorems is that both periodically oscillating functional $F$ and the corresponding effective functional $\bar F$ admit interior $C^{2,\gamma}$ estimates when the coefficients are fixed. More precisely, $F,\bar F \in S_2$ in the sense of Definition \ref{definition:class}. 

This condition becomes straightforward, when $F$ is a concave functional. First, we have $F\in S_2$ by the Evans-Krylov theory \cite[Theorem 6.1]{CC}. Next, according to \cite[Lemma 3.2]{E}, if $F$ is concave, then so is $\bar F$, proving $\bar F\in S_2$ by the same theory again. 

Henceforth, we shall find a class of non-concave functionals $F$ that both $F$ and $\bar F$ belong to class $S_2$. This will imply that the class of periodically oscillating functionals $F$ that verify the assumptions of Theorem \ref{theorem:int-C11} and Theorem \ref{theorem:bdry-C11} strictly wider than the class of concave functionals. 

Let us being with a lower bound for the Hessian of the interior corrector. 

\begin{lemma}\label{lemma:wF-C11}
Let $F\in S_2(\lambda,\Lambda,\bar\kappa,\bar\gamma)\cap R_0(\kappa,\gamma)$ be a functional on $\cS^n\times \R^n$ that is periodic in the second argument. Then there is some constant $L>1$, depending only on $n$, $\lambda$, $\Lambda$, $\bar\kappa$, $\bar\gamma$, $\kappa$ and $\gamma$, such that for any $M\in\cS^n$, one has
\begin{equation}\label{eq:wF-C11}
\min_{\R^n} |D_y^2 w_F(M,\cdot) + M| \geq \frac{|M|}{L},\quad\text{if }|M| > L. 
\end{equation} 
\end{lemma} 

\begin{proof}
Let us fix $M\in\cS^n$, and denote by $w$ the interior corrector $w_F(M,\cdot)$. According to \cite{E}, $w$ is the limit function of sequence $\{w^\delta - w^\delta(0)\}_{\delta>0}$, where $w^\delta$ is the unique, periodic viscosity solution to  
\begin{equation*}
F(D_y^2 w^\delta + M,y) -\delta w^\delta = 0 \quad\text{in }\R^n.
\end{equation*}
In particular, constant functions $\delta^{-1} \min_{\R^n} F(M,\cdot)$ and $\delta^{-1} \max_{\R^n} F(M,\cdot)$ are a periodic viscosity subsolution and respectively supersolution to the penalized problem above. Owing to this fact, one can deduce that
\begin{equation}\label{eq:Fb-re} 
\min_{\R^n} F(M,\cdot) \leq \bar{F}(M)\leq \max_{\R^n} F(M,\cdot).
\end{equation} 
For this reason, one can also derive a sharper interior $C^{2,\gamma}$ estimate
\begin{equation}\label{eq:w-C2a-re}
\norm{w}_{C^{2,\gamma}(\R^n)} \leq C_1 \norm{F(M,\cdot)}_{L^\infty(\R^n)},
\end{equation} 
compared to \eqref{eq:wF-C2a}, where $C_1>0$ depends only on $n$, $\lambda$, $\Lambda$, $\bar\kappa$, $\bar\gamma$, $\kappa$ and $\gamma$. 

Let $\chi>1$ be a given number, fix $L>1$ by a sufficiently large number to be determined at the end of the proof, and suppose that 
\begin{equation}\label{eq:wF-C11-assmp}
|M| > L\chi.
\end{equation} 
Let us first consider the case where 
\begin{equation}\label{eq:F-G-case1}
\norm{F(M,\cdot)}_{L^\infty(\R^n)} < \frac{(L-1)}{C_1} \chi. 
\end{equation}
Then it follows from \eqref{eq:w-C2a-re} that 
\begin{equation*}
|D^2 w + M| \geq |M| - |D^2 w| \geq L\chi - C_1 \norm{F(M,\cdot)}_{L^\infty(\R^n)} > \chi.
\end{equation*} 
This proves \eqref{eq:wF-C11} under the additional assumption \eqref{eq:F-G-case1}. 

Next, let us consider the other case where
\begin{equation}\label{eq:F-G-case2}
\norm{F(M,\cdot)}_{L^\infty(\R^n)} \geq \frac{(L-1)}{C_1}\chi.
\end{equation}
Due to the assumption $F\in R_0(\kappa,\gamma)$, we have $\osc_{\R^n} F(M,\cdot) \leq \kappa \sqrt{n}$. Hence, taking $L$ sufficiently large such that $(L-1) \chi > C_1\kappa \sqrt{n}$, then owing to \eqref{eq:F-G-case2}, we may assume loss of generality that 
\begin{equation*}
\min_{\R^n} F(M,\cdot) \geq \frac{L-1}{C_1} \chi - \kappa \sqrt{n} > 0.
\end{equation*}
Therefore, it follows from \eqref{eq:cell}, \eqref{eq:Fb-re} and the ellipticity assumption $F \in S_0(\lambda,\Lambda)$ that 
\begin{equation}\label{eq:F-G-more}
\frac{L-1}{C_1} \chi - \kappa \sqrt{n} \leq F(D^2 w + M,y) \leq \Lambda |D^2 w + M|\quad\text{in }\R^n. 
\end{equation}
Finally, we choose $L$ large enough so as to satisfy 
\begin{equation}\label{eq:L-cond}
\frac{L - 1}{C_1} - \kappa\sqrt{n} > \Lambda;
\end{equation} 
note that $L$ depends only on $\Lambda$, $\kappa$ and $C_1$, hence on $n$, $\lambda$, $\Lambda$, $\kappa$ and $\gamma$ only. Then one may verify from \eqref{eq:F-G-more} as well as the assumption $\chi > 1$ that 
\begin{equation*}
|D^2 w + M| > \chi\quad\text{in }\R^n,
\end{equation*}
again proving \eqref{eq:wF-C11}. Hence, we have verified that \eqref{eq:wF-C11} holds under the other assumption \eqref{eq:F-G-case2}, from which we conclude that it holds in general. This finishes the proof. 
\end{proof} 

Next, we present some monotone property of effective functionals. 

\begin{lemma}\label{lemma:mono}
Let $F_1$ and $F_2$ be uniformly elliptic, periodic and continuous functionals on $\cS^n\times\R^n$. Then 
\begin{equation}\label{eq:mono}
\overline{\min\{F_1, F_2\}} \leq \min\{\bar F_1, \bar F_2\}\quad\text{on }\cS^n. 
\end{equation}  
\end{lemma} 

\begin{remark}\label{remark:mono}
In general, we do not have equality in \eqref{eq:mono}, even if $F_1$ and $F_2$ are linear functionals. For example, if $n=1$, and $F_1(M,y) = ( 2 + \cos(2\pi y))M$, $F_2(M,y) = (2 + \sin(2\pi y))M$, then 
\begin{equation*}
\bar F_1(M) = \left( \int_0^1 \frac{dy}{2 + \cos(2\pi y)} \right)^{-1} M = \left(\int_0^1 \frac{dy}{2 + \sin(2\pi y)}\right)^{-1} M = \bar F_2(M), 
\end{equation*} 
for any $M\in \cS^1 = \R$, and therefore,
\begin{equation*}
\overline{\min\{ F_1,F_2\}}(M) = \left( \int_0^1 \frac{dy}{2 + \min\{\cos(2\pi y),\sin(2\pi y)\}} \right)^{-1} M < \bar F_1(M) = \bar F_2(M), 
\end{equation*}
for any $M\neq 0$. 
\end{remark}

\begin{proof}[Proof of Lemma \ref{lemma:mono}]
Set $F_0 = \min\{F_1,F_2\}$. Suppose towards a contradiction that \eqref{eq:mono} fails at some $M\in\cS^n$. For each $i\in\{0,1,2\}$, denote by $w_i$ and $\gamma_i$ the interior corrector $w_{F_i}(M,\cdot)$ and respectively the value $\bar F_i(M)$ of the effective functional corresponding to $F_i$, in the sense of Definition \ref{definition:eff}. Then for each $i\in\{1,2\}$, 
\begin{equation}\label{eq:mono-fail}
F_i (M + D^2 w_0,y) \geq F_0(M + D^2 w_0,y) = \gamma_0 > \gamma_i\quad\text{in }\R^n,
\end{equation}
in the viscosity sense. As $w_i$ being a periodic viscosity solution to $F_i (M + D^2 w,y) = \gamma_i$ in $\R^n$, we deduce that 
\begin{equation*}
P^+(D^2 (w_0 - w_i)) > 0\quad\text{in }\R^n, 
\end{equation*} 
in the viscosity sense, where $P^+$ is the Pucci maximal operator associated with the ellipticity bounds for both $F_0$ and $F_i$. However, as both $w_0$ and $w_i$ being periodic, $w_0 - w_i$ attains a maximum at some $y_i \in Q_1$, where $Q_1$ is the unique periodic cube, so 
\begin{equation*}
P^+ (D^2 (w_0 - w_i)(y_i)) \leq 0,
\end{equation*} 
in the viscosity sense, a contradiction. 
\end{proof}

Our strategy is as follows. We shall use the Evans-Krylov theory developed by Cabr\'e and Caffarelli \cite{CC2} for the class of non-concave functionals $F$ given by the minimum of a concave and a convex functional, say $F^\cap$ and respectively $F^\cup$. The advantage of this class is that those functionals $F$ admit the interior $C^{2,\bar\gamma}$ estimates with constant $\bar\kappa$ for certain $\bar\kappa>1$ and $\bar\gamma\in(0,1)$ that depend only on the dimension and the ellipticity constants, whence the estimates do not change under translation and scaling, i.e., $F \in S_2(\lambda,\Lambda,\bar\kappa,\bar\gamma)$ from \cite[Theorem 1.1]{CC2}; recall that class $S_2(\lambda,\Lambda,\bar\kappa,\bar\gamma)$ requires that not only $F$ but also all of its the translated versions, $(M,y)\mapsto (F(M+N,y) - F(M,y))$, satisfy the interior $C^{2,\bar\gamma}$ estimates with the same constant $\bar \kappa$.  

The question is if $\bar F \in S_2(\lambda,\Lambda,\bar c,\bar\alpha)$ holds for some $\bar c>1$ and $\bar\alpha\in(0,1)$. Here, we shall impose some additional conditions on $F^\cap$ and $F^\cup$, whose minimum produces $F$, such that $\bar F = \min\{\bar F^\cap,\bar F^\cup\}$. As a byproduct, $\bar F\in S_2(\lambda,\Lambda,\bar\kappa,\bar\gamma)$ with $\bar\kappa$ and $\bar\gamma$ as above. 

\begin{remark}\label{remark:key}
Let us stress that the condition $\bar F = \min\{\bar F^\cap,\bar F^\cup\}$ is sufficient to have $\bar F \in S_2$, but by no means necessary. A good example is shown in Remark \ref{remark:mono}: with the specific choice of the linear functionals $F^\cap$ and $F^\cup$ there, we have $\bar F < \min\{\bar F^\cap,\bar F^\cup\}$ on $\cS^n\setminus\{0\}$, while $F$ as the minimum of two linear functionals is concave, and so is $\bar F$, from which we obtain $\bar F \in S_2$. 
\end{remark} 

\begin{lemma}\label{lemma:key}
Let $F^\cap$ and $F^\cup$ be uniformly elliptic, periodic and continuous functionals on $\cS^n\times\R^n$, with ellipticity constants $\lambda$ and $\Lambda$, such that $F^\cap$ is concave, while $F^\cup$ is convex in the first argument, and that $\{F^\cap, F^\cup\}\subset R_0(\kappa,\gamma)$ for some $\kappa>0$ and $\gamma\in(0,1)$. There exists some $L>1$, depending only on $n$, $\lambda$, $\Lambda$, $\kappa$ and $\gamma$, such that if $F^\cap$ and $F^\cup$ also satisfy, for some $R>0$,
\begin{equation}\label{eq:key-const}
 \osc_{\R^n} F^\cap (M,\cdot) = \osc_{\R^n} F^\cup (M,\cdot) = 0,\quad\text{whenever } |M|\leq LR\text{ and}
\end{equation} 
\begin{equation}\label{eq:key-large}
\min_{\R^n} \left( F^\cup (M,\cdot) - F^\cap (M,\cdot) - \kappa n^{\frac{\gamma}{2}} |M|\right) \geq 0,\quad\text{whenever } |M|\geq R,
\end{equation} 
then 
\begin{equation}\label{eq:key} 
\overline{\min\{ F^\cap, F^\cup\}}= \min\{ \bar F^\cap , \bar F^\cup\}\quad\text{on }\cS^n. 
\end{equation}   
\end{lemma} 

\begin{proof}
Define $F = \min\{F^\cap, F^\cup\}$ on $\cS^n\times\R^n$. 
The first condition \eqref{eq:key-const} implies that $\osc_{\R^n} F(M,\cdot) = 0$ for all $M\in\cS^n$ with $|M|\leq LR$. Hence, we have $\bar F(M) = F(M,\cdot)$, $\bar F^\cap (M) = F^\cap (M,\cdot)$ and $\bar F^\cup (M) = G^\cup(M,\cdot)$ on $\R^n$. In particular, \eqref{eq:key} holds for all $M\in\cS^n$ with $|M| \leq LR$. 

Hence, we are left with proving \eqref{eq:key} for $M\in\cS^n$ with $|M| > LR$, with $L$ to be determined as in Lemma \ref{lemma:wF-C11}. To be more precise, we choose $L_1$ as follows. Since $\lambda$ and $\Lambda$ are the ellipticity bounds for both $F^\cap$ and $F^\cup$, the Evans-Krylov theory \cite[Theorem 6.1]{CC} for concave functionals implies that $\{F^\cap - F^\cap(0,\cdot), F^\cup - F^\cup(0,\cdot)\} \subset S_2(\lambda,\Lambda,\bar\kappa,\bar\gamma)$ for some $\bar\kappa>1$ and $\bar\gamma\in(0,1)$, depending only on $n$, $\lambda$ and $\Lambda$. With such $\bar\kappa$ and $\bar \gamma$, let us select $L>1$ by the large constant as in Lemma \ref{lemma:wF-C11} such that \eqref{eq:wF-C11} holds for any $G \in S_2(\lambda,\Lambda,\bar\kappa,\bar\gamma)\cap R_0(\kappa,\gamma)$, in particular for any $G\in \{F-F(0,\cdot), F^\cap - F^\cap(0,\cdot), F^\cup - F^\cup(0,\cdot)\}$. Due to the universality of $\bar\kappa$ and $\bar\gamma$, $L$ depends only on $n$, $\lambda$, $\Lambda$, $\kappa$ and $\gamma$. 

Fix any $M\in\cS^n$ such that $|M| \geq LR$. For the sake of simplicity, denote by $w$, $w_\cap$ and $w_\cup$ the interior corrector $w_F(M,\cdot)$, $w_{F^\cap}(M,\cdot)$ and respectively $w_{F^\cup}(M,\cdot)$ in the sense of Definition \ref{definition:eff}. Since Lemma \ref{lemma:wF-C11} applies to $F$, we observe from \eqref{eq:wF-C11} that 
\begin{equation*}
\min_{\R^n} |D_y^2 w + M| \geq \frac{|M|}{L} > R, 
\end{equation*} 
so the assumption \eqref{eq:key-large} implies that $F^\cap(D^2 w(y) + M,y) \leq F^\cup(D^2 w (y)+ M,y) - \kappa n^{\gamma/2}|M| < F^\cup(D^2 w(y) + M,y)$ for all $y\in\R^n$. Thus, $F(D^2 w(y) + M,y) = F^\cap(D^2 w(y) + M,y)$ for all $y\in\R^n$, which in turn yields that 
\begin{equation*}
F^\cap (D^2 w + M,y) = \bar F(M)\quad\text{in }\R^n. 
\end{equation*} 
As $w$ being a periodic function and $\bar F(M)$ being a constant, it follows from the uniqueness of the effective functional that 
\begin{equation}\label{eq:key1}
\bar F (M) = \bar F^\cap (M).
\end{equation} 

To this end, we claim that $\bar F^\cap (M) \leq \bar F^\cup(M)$ if $|M|\geq LR$.  Note that $w_\cap$ and $w_\cup$ are periodic functions on $\R^n$, and that $w_\cap, w_\cup \in C^{2,\gamma}(\R^n)$; these are ensured by the fact that $\{F^\cap - F^\cap(0,\cdot), F^\cup - F^\cup(0,\cdot)\}\subset S_2(\lambda,\Lambda,\bar\kappa,\bar\gamma)\cap R_0(\kappa,\gamma)$. Hence, we can find some $y_\cap,y_\cup \in \bar Q_1 = [-\frac{1}{2},\frac{1}{2}]^n$ be such that $|D^2 w_\cap (y_\cap)| = |D^2 w_\cup(y_\cup)| = 0$. In view of the assumption \eqref{eq:key-large} and $F^\cup\in R_0(\kappa,\gamma)$, as well as the cell problems that $w_\cap$ and $w_\cup$ solve, we derive that 
\begin{equation*}
\bar F^\cap (M) = F^\cap ( M, y_\cap) \leq F^\cup (M, y_\cap) - \kappa n^{\gamma/2} |M| \leq F^\cup (M, y_\cup) = F^\cup (M).
\end{equation*} 
This together with \eqref{eq:key1} finishes the proof of \eqref{eq:key} for $M\in\cS^n$ with $|M|\geq LR$. 
\end{proof} 

Consequently, we obtain a class of non-concave, periodic functionals whose effective functionals admit interior $C^{2,\alpha}$ estimate. In particular, the nonlinear homogenization problem with this class of functionals admit uniform interior $C^{1,1}$ estimates. 

\begin{proposition}\label{proposition:ex}
Let $F^\cap$ and $F^\cup$ be as in Lemma \ref{lemma:key}, and define $F = \min\{F^\cap, F^\cup\}$. Suppose that $F(0,\cdot) = 0$ on $\R^n$. Then $F\in S_2(\lambda,\Lambda,\bar\kappa,\bar\gamma)\cap R_0(\kappa,\gamma)$, and $\bar F \in S_2(\lambda,\Lambda,\bar\kappa,\bar\gamma)$, where $\kappa>0$ and $\gamma\in(0,1)$ are as in Lemma \ref{lemma:key}, while $\bar\kappa>1$ and $\bar\gamma\in(0,1)$ are some constants depending only on $n$, $\lambda$ and $\Lambda$. Moreover, if $f\in C^\alpha(B_1)$ and $u^\e\in C(B_1)$ is a viscosity solution of 
\begin{equation*}
F \left( D^2 u^\e,\frac{x}{\e}\right) = f \quad\text{in }B_1,
\end{equation*} 
then $u^\e \in C^{1,1}(B_{1/2})$ and 
\begin{equation*}
\norm{u^\e}_{C^{1,1}(B_{1/2})} \leq C \left( \norm{u^\e}_{L^\infty(B_1)} + \norm{f}_{C^\alpha(B_1)} \right), 
\end{equation*} 
where $C>0$ depends only on $n$, $\lambda$, $\Lambda$, $\bar\kappa$, $\bar\gamma$, $\kappa$, $\gamma$ and $\alpha$. 
\end{proposition}

\begin{proof}
The fact that $F\in S_2(\lambda,\Lambda,\bar\kappa,\bar\gamma)\cap R_0(\kappa,\gamma)$ follows immediately from the assumption on $F^\cap$ and $F^\cup$ in Lemma \ref{lemma:key} as well as the discussion above the statement of the lemma. On the other hand, since the effective functional to a concave/convex functional is also concave/convex \cite[Lemma 3.2]{E}, it follows from \eqref{eq:key} that $\bar F$ is the minimum of a concave and a convex, homogeneous functional. Thus, we infer from \cite[Theorem 1.1]{CC2} that $\bar F \in S_2(\lambda,\Lambda,\bar\kappa,\bar\gamma)$. Therefore, $F$ and $\bar F$ verify the assumptions of Theorem \ref{theorem:int-C1a}, from which the uniform interior $C^{1,1}$ estimate of $u^\e$ follows. 
\end{proof} 

\begin{remark}\label{remark:ex}
We do not claim a uniform boundary $C^{1,1}$ estimate for $F$ as in Proposition \ref{proposition:ex}. Note that $F$, as the minimum of a concave and a convex functional, is Lipschitz in the matrix variable in general. Hence, $F\in R_0\setminus R_1$, while Theorem \ref{theorem:bdry-C11} requires $F$ to belong to $R_1$; recall from Definition \ref{definition:class} that $R_1$ consists of functionals whose derivatives in the matrix variable are Lipschitz.
\end{remark}


\end{document}